\documentclass[12pt,a4paper]{article}

\usepackage{amssymb}
\usepackage{theorem}
\usepackage[dvips]{epsfig}
\usepackage{latexsym}
\usepackage{exscale}
\usepackage[latin1]{inputenc}
\newtheorem{thm}{Theorem}[section]
\newtheorem{lemma}[thm]{Lemma}
\newtheorem{prop}[thm]{Proposition}
\newtheorem{fact}{Fact}[section]

\newtheorem{cor}[thm]{Corollary}
\newcommand{\pf}{\noindent {\it Proof:\ }}
\newcommand{\remark}{\noindent {\bf Remark:\ }}

\def\be#1\ee{\begin{equation}#1\end{equation}}
\newcommand{\bea}{\begin{eqnarray}}
\newcommand{\eea}{\end{eqnarray}}
\newcommand{\beaa}{\begin{eqnarray*}}
\newcommand{\eeaa}{\end{eqnarray*}}
\newcommand{\bei}{\begin{itemize}}
\newcommand{\eei}{\end{itemize}}
\newcommand{\bee}{\begin{enumerate}}
\newcommand{\eee}{\end{enumerate}}

\def\norm#1{\left\|#1\right\|}             
\def\abs#1{\left\vert #1 \right\vert}      

\def\set#1{\left\{#1\right\}}

\def\P{{\mathbb{P}}}
\def\pr#1{\P\left(#1\right)}

\def\e{{\mathbb E}}

\def\n{{\mathbb N}}

\def\ee{\mathrm{e}}
\def\d{\, \mathrm{d}}

\def\diam{{\mathrm{diam}}}

\def\E{{\mathbb{E}}\,}
\def\R{{\mathbb{R}}}
\def\N{\mathbb{N}}
\def\H{\mathcal{H}}
\def\Z{\mathbb{Z}}

\def\on{{\mathbf 1}}
\def\dl{{\Delta}}
\def\ph{{\varphi}}

\newcommand{\eps}{\varepsilon}

\def\ee{\mathrm{e}}      
\def\d{\,\mathrm{d}}     

\begin{document}
\title{Small Deviations of Gaussian Random Fields
in $L_q$--Spaces}

\author{Mikhail Lifshits\and Werner Linde \and Zhan Shi
   } \date{} \maketitle

\bigskip

\begin{abstract}
\noindent

We investigate small deviation properties of Gaussian random
fields
in the space $L_q(\R^N,\mu)$ where $\mu$ is an arbitrary
finite compactly
supported Borel measure. Of special interest are hereby ``thin"
measures $\mu$, i.e., those which are singular with respect to
the $N$--dimensional Lebesgue measure; the so--called self--similar measures providing
a class of typical examples.

For a large class of random fields
(including, among others, fractional Brownian motions),
we describe the behavior of small deviation probabilities
via numerical characteristics of $\mu$, called
mixed entropy, characterizing size and regularity of $\mu$.

For the particularly interesting
case of self--similar measures $\mu$, the asymptotic behavior of
the mixed entropy is evaluated explicitly. As a consequence, we get the
asymptotic of the small deviation for $N$--parameter fractional Brownian motions
with respect to $L_q(\R^N,\mu)$--norms.

While the upper estimates for the small
deviation probabilities are proved by purely probabilistic methods,
the lower bounds are established by analytic tools concerning
Kolmogorov and entropy numbers of H\"older operators.

\end{abstract}
\bigskip
\bigskip
\bigskip
\bigskip
\bigskip

\vfill

\noindent

{\bf Key words:}\  random fields, Gaussian processes,
 fractional Brownian motion,
fractal measures, self--similar measures, small deviations,
Kolmogorov numbers, metric entropy, H\"older operators.

\baselineskip=6.0mm

\newpage

\section{Introduction}


The aim of the present paper is the investigation of the small deviation
behavior of Gaussian random fields in the $L_q$--norm
taken with respect to a rather arbitrary measure on $\R^N$.
Namely, for a Gaussian random field $(X(t),t\in \R^N)$, for a measure
$\mu$ on $\R^N$, and for any $q\in [1,\infty)$ we are interested
in the behavior of the small deviation function
\be
\label{sdf}
\ph_{q,\mu}(\eps) :=
- \log \pr{\int_{\R^N}\abs{X(t)}^q \d \mu(t)<\eps^q}\,,
\ee
as $\eps\to 0$ in terms of certain quantitative properties of the
underlying measure
$\mu$.
Let us illustrate this with an example. As a consequence of our estimates, we get
the following corollary for the $N$--parameter fractional
Brownian motion $W_H= (W_H(t), t\in\R^N)$. For the exact meaning of the theorem, see Section \ref{s:ss}.

\begin{thm} \label{th_exammple}
 Let $T\subset \R^N$ be a compact self--similar set of Hausdorff dimension $D>0$ and
let $\mu$ be the $D$-dimensional Hausdorff measure on $T$.
Then
$$
- \log \pr{\int_{T}\abs{W_H(t)}^q \d \mu(t)<\eps^q} \approx \eps^{-D/H}\;.
$$
\end{thm}

General small deviation problems attracted much attention during the
last years due to their deep relations to various mathematical
topics such as operator theory, quantization, strong limit laws in
statistics, etc, see the surveys \cite{LS,Lif}. A more specific
motivation for this work comes from \cite{LLS}, where the
one--parameter case $N=1$ was considered for fractional Brownian
motions and Riemann--Liouville processes.
\medskip

Before stating our main multi--parameter results, let us recall a basic
theorem from \cite{LLS}, thus giving a clear idea of the entropy
approach to small deviations in more general $L_q$--norms.

Recall that the (one--parameter) fractional Brownian motion (fBm) $W_H$
with Hurst index $H\in (0,1)$
is a
centered Gaussian process on $\R$ with a.s.~continuous paths
and covariance
$$
\E W_H(t)\,W_H(s) = \frac{1}{2}\set{t^{2 H}+ s^{2 H}-\abs{t-s}^{2 H}}\;,
\quad t,s\in \R .
$$

We write $f\sim g$ if
  $\lim_{\eps\to 0} \frac{f(\eps)}{g(\eps)}=1$ while
$f\preceq  g$ (or $g \succeq f$) means that
   $\limsup_{\eps\to 0} \frac{f(\eps)}{g(\eps)}<\infty$. Finally,
$f\approx g$ says that $f\succeq g$ as well as $g\succeq f$.

If  $\mu=\lambda_1$, the restriction of the Lebesgue
measure to $[0,1]$, then for $W_H$
the behavior
of $\ph_{q,\mu}(\eps)$ is well--known, namely,
$\ph_{q,\mu}(\eps)\sim c_{q,H}\, \eps^{-1/H}$, as $\eps\to 0$.
The exact value of the finite and positive constant $c_{q,H}$
is known only in few cases; sometimes a variational representation
for $c_{q,H}$ is available. See more details in \cite{LS} and
\cite{LSi}.

If $\mu$ is absolutely continuous with respect to $\lambda_1$,
the behavior of $\ph_{q,\mu}(\eps)$ was investigated in \cite{Li1},
\cite{LLMem} and \cite{LLTAMS}. Under mild assumptions, the order
 $\eps^{-1/H}$ remains unchanged, only
an extra factor depending on the density of $\mu$ (with respect
to $\lambda_1$) appears.
The situation is completely different for measures $\mu$ being
singular  to $\lambda_1$. This question was recently
investigated in \cite{LJAT} for $q=\infty$ (here only the size of
the support of $\mu$ is of importance) and in \cite{Na} for
self--similar measures  and $q=2$.
When passing from $q=\infty$ to a finite $q$, the problem becomes more
involved because in the latter case the distribution of the mass
of $\mu$ becomes important.
Consequently, one has to introduce some kind of entropy of $\mu$
taking into account the size of its support as well as the
distribution of the mass on $[0,1]$. This is done in the following way.

Let $\mu$ be a continuous measure on $[0,1]$, let $H>0$ and $q\in[1,\infty)$.
We define a number $r>0$ by
\be
\label{r1}
1/r := H + 1/q \; .
\ee
Given an
interval $\dl\subseteq [0,1]$, we denote
\be \label{J}
J_{\mu}^{(H,q)}(\dl):=\abs{\dl}^H\cdot\mu(\dl)^{1/q}
\ee
and set
\be
\label{sm}
\sigma_\mu^{(H,q)}(n) :=
\inf\set{\left(\sum_{j=1}^n J_\mu(\dl_j)^r\right)^{1/r}
 :\; [0,1]\subseteq \bigcup_{j=1}^n \dl_j}\; ,
\ee
where the $\dl_j$'s are supposed to be intervals on the real line.
The sequence $\sigma_\mu^{(H,q)}(n)$ may be viewed as some kind of outer
{\it mixed} entropy of $\mu$. Here ``mixed" means that we take
into account  the measure as well as the length of an interval.

The main result of \cite{LLS} shows a very tight relation between
the behavior of $\sigma_\mu^{(H,q)}(n)$, as $n\to\infty$, and of the small
deviation function (\ref{sdf}). More precisely, the following is true.

\begin{thm} \label{lls}
Let $\mu$ be a finite continuous measure on $[0,1]$ and let $W_H$
be a fBm of index $H>0$. For
$q\in[1,\infty)$, define $\sigma_\mu^{(H,q)}(n)$ as in $(\ref{sm})$.
\bee
\item[\rm (a)]
If
$$
 \sigma_\mu^{(H,q)}(n)\succeq n^{-\nu}\,(\log n)^\beta
$$
for certain $\nu\ge 0$ and $\beta\in\R$, then
$$
 - \log\P\{\norm{W_H}_{L_q([0,1],\mu)} < \eps\}
 \succeq
 \eps^{-1/(H+\nu)}\cdot\log(1/\eps)^{\beta/(H+\nu)}\;.
$$

\item[\rm (b)]
On the other hand, if
$$
\sigma_\mu^{(H,q)}(n)\preceq n^{-\nu}\,(\log n)^\beta
$$
then
\be
\label{low}
- \log\P\{\norm{W_H}_{L_q([0,1],\mu)}<\eps\}
 \preceq
 \eps^{-1/(H+\nu)}\cdot\log(1/\eps)^{\beta/(H+\nu)}\;.
\ee
\eee
\end{thm}
\bigskip

Remarkably, there is another quantity, a kind of inner mixed
entropy, equivalent to $\sigma_\mu^{(H,q)}(n)$ in the one--parameter case.
This one is defined as follows. Given $\mu$ as before, for each
$n\in\N$ we set
\be \label{dem}
 \delta_\mu^{(H,q)}(n):=
 \sup\set{\delta>0 : \exists\, \Delta_1,\ldots,\Delta_n\subset [0,1],\; J_\mu^{(H,q)}(\Delta_i)\ge \delta}
\ee
where the $\Delta_i$ are supposed to possess disjoint interiors.

It is shown in \cite{LLS} that  $\sigma_\mu^{(H,q)}(n)$ and
$ n^{1/r}\delta_\mu^{(H,q)}(n)$ are, in a sense, equivalent as
$n\to\infty$, namely, it is proved that
for each integer $n\ge 1$, we have
\be  \label{q1}
  \sigma_\mu^{(H,q)}(2 n+1)\le (2 n+1)^{1/r}\, \delta_\mu^{(H,q)}(n)
  \quad\mbox{and}\quad n^{1/r}\,\delta_\mu^{(H,q)}(2n)\le \sigma_\mu^{(H,q)}(n)\;.
\ee
Therefore,  Theorem \ref{lls} can be immediately restated
in terms of $\delta_\mu^{(H,q)}(n)$. For example, $\delta_\mu^{(H,q)}(n) \approx
n^{-(1/q+1/a)}(\log n)^{\beta}$ with $a\le 1/H$ is equivalent to
$\sigma_\mu^{(H,q)}(n)\approx n^{-(1/a-H)}(\log n)^{a\beta}$
and thus to
$$
- \log\P\{\norm{W_H}_{L_q([0,1],\mu)}<\eps\}
 \approx
 \eps^{-a}\cdot (\log(1/\eps))^{a\beta}\;.
$$
Notice that in the case of measures on $[0,1]$
the restriction $a\le 1/H$ is natural; it is attained for the
Lebesgue measure.

\medskip

To our great deception, we did not find in the literature the
notions of outer and inner mixed entropy as defined above, although
some similar objects do exist: cf.~the notion of weighted Hausdorff measures
investigated in \cite{Ma}, pp.~117--120,
or $C$-structures associated with metrics and measures in
Pesin \cite{Pe}, p.~49, or
multifractal generalizations of Hausdorff measures and packing measures in Olsen \cite{Ol}.
Yet the quantitative properties of Hausdorff dimension and entropy of
a set seem to have  almost nothing in common (think of any countable set ---
its Hausdorff dimension is zero while the entropy
properties can be quite non--trivial).

\section {Main results in the multi--parameter case}
\setcounter{equation}{0}

Although in this article we are mostly interested in the behavior of the
$N$--parameter fractional Brownian motion, an essential part of our estimates
is valid for much more general processes.
For example, to prove lower estimates for $\ph_{q,\mu}(\eps)$
 we only need a certain non-degeneracy of interpolation
errors (often called ``non-determinism") while for upper estimates of $\ph_{q,\mu}(\eps)$ some H\"older
type inequality suffices. Therefore, we start from the
general setting.
Let $X :=(X(t), t\in T)$ be a centered measurable Gaussian process on a metric
space $(T,\rho)$. Here we endow $T$ with the $\sigma$--algebra of Borel sets.
For $t\in T$, any $A\subset T$ and  $\tau>0$
we set
\be
\label{vt}
    v(t ,\tau):=\left( \mathrm {Var}
     \left[ X(t) -
     \E  \left( X(t) \, | \, X(s), \, \rho(s,t) \ge \tau
        \right)
     \right]\right)^{1/2}
\ee
and
\be
\label{vt2}
      v(A,\tau):=  \inf_{t\in A} v(t,\tau)\;.
\ee
We suppose that $\mu$ is a finite Borel measure on $T$ and that
$A_1$, $\dots$, $A_n$ are disjoint measurable subsets of $T$. Let
$$
   v_1 :=  \inf_{t\in A_1} \left(\mathrm {Var} X(t)\right)^{1/2}\cdot \mu(A_1)^{1/q},
$$
and
\begin{equation}
    v_i := v \big(A_i , \tau_i\big)\cdot \mu(A_i)^{1/q}, \qquad
    2\le i\le n\,,
    \label{Jbar}
\end{equation}
where $\tau_i :=\mathrm{dist} (A_i, \,
    \bigcup_{k=1}^{i-1} A_k )$. We set
\be
\label{wmu}
V_\mu := V_\mu(A_1, \dots, A_n) := \min_{1\le i\le n}
v_i.
\ee
Finally, given $n\in\N$ we define some kind of weighted inner entropy by
\be
\label{tilded}
\delta_\mu(n):= \sup\set{\delta>0 : \exists\;\textrm{disjoint}\;A_1,\ldots,A_n\subset T,\;
V_\mu(A_1, \dots, A_n)\ge\delta}\;.
\ee

\noindent
In this quite general setting we shall prove the following.

\medskip

\begin{thm}
\label{up1}
Let $\mu$ be a finite measure on $T$ and let $X$
be a measurable centered Gaussian random field on $T$. Let
$q\in[1,\infty)$ and define $\delta_\mu(n)$ as in $(\ref{tilded})$.
If
$$
\delta_\mu(n)\succeq n^{-1/q-1/a}\,(\log n)^\beta
$$
for certain $a>0$ and $\beta\in\R$, then
$$ 
 - \log\P\{\norm{X}_{L_q(T,\mu)} < \eps\}
 \succeq
 \eps^{-a}\cdot\log(1/\eps)^{a\beta}\;.
$$
\end{thm}
\bigskip

For the case $T\subset \R^N$ ($N\ge 1$) with the metric $\rho$ generated by the
Euclidean distance, i.e., $\rho(t,s)=\abs{t-s}$, $t,s\in\R^N$, we give a slightly weaker upper
bound for the small deviation probabilities. This bound, however, has the
advantage of using simpler geometric characteristics. In particular,
we do not need to care about the distances between the sets.
If $A\subset T$ is measurable, we set
$$ 
\overline v (A) := v\left(A, \tau_A \right)
$$
where $\tau_A:= \mathrm{diam} (A)/ (2\sqrt{N})$, i.e.,
$$
\overline v (A)=\inf_{t\in A}\left( \mathrm {Var}
     \left[ X(t) -
     \E  \left( X(t) \, | \, X(s), \,s\in T,\;2\sqrt N \abs{t-s} \ge \mathrm{diam} (A)
        \right)
     \right]\right)^{1/2}\;.
$$
Given cubes $Q_1,\ldots,Q_n$ in $T$ with disjoint interiors, similarly as in (\ref{wmu}),
we define the quantity
$$
\overline V_\mu= \overline V_\mu(Q_1,\ldots,Q_n):= \inf_{1\le i\le n} \overline v(Q_i)\cdot
\mu(Q_i)^{1/q}
$$
and as in (\ref{tilded}) we set
\be \label{demd}
\overline \delta_\mu(n):= \sup\set{\delta>0 : \exists\,Q_1,\ldots,Q_n\subset T,\;
\overline V_\mu(Q_1, \dots, Q_n)\ge\delta}
\ee
where the cubes $Q_i$ are supposed to possess disjoint interiors.
\medskip

We shall prove the following.
\begin{thm} \label{t:upd}
Let $\mu$ be a finite measure on $T\subset \R^N$ and let $X$
be a centered Gaussian random field on $T$. For
$q\in[1,\infty)$ and $\overline \delta_\mu(n)$ defined as in $(\ref{demd})$, if
$$
\overline \delta_\mu(n)\succeq n^{-1/q-1/a}\,(\log n)^\beta
$$
for certain $a>0$ and $\beta\in\R$, then
$$ 
 - \log\P\{\norm{X}_{L_q(T,\mu)} < \eps\}
 \succeq
 \eps^{-a}\cdot\log(1/\eps)^{a\beta}\;.
$$
\end{thm}
\bigskip

Finally we apply our results to the $N$--parameter fractional Brownian motion,
i.e., to the
real--valued centered Gaussian random field $W_H:= (W_H(t), \, t\in \R^N)$
with covariance
$$
\e\left[ W_H(s) W_H(t) \right] = {1\over 2} \left(  | s|^{2H} + |
t|^{2H} - |t-s|^{2H} \right) , \qquad (s,t) \in \R^N \times \R^N,
$$

\noindent where $H\in (0,1)$ is the Hurst index.

It is known that $W_H$ satisfies (see \cite{pitt} and \cite{xiao} for further information about
processes satisfying similar conditions)
\begin{equation}
\label{LND}
    \mathrm{Var} \left[ W_H(t)- \E\left( W_H(t) \, | \, W_H(s), \, |s-t|
\ge \tau \right) \right]
\ge c\ \tau^{2H}, \qquad t\in \R^N,
\end{equation}
for all $0\le \tau\le |t|$.
Thus, in view of (\ref{LND}) it is rather natural to adjust the definition of
$\overline \delta_\mu$ as follows.
Namely, as in (\ref{J}) for $N=1$,   we set
\be
\label{Jmu}
J^{(H,q)}_\mu(A) := \left(\mathrm{diam}(A)\right)^{H} \mu(A)^{1/q}
\ee
for any measurable subset $A\subset\R^N$.
Then the multi--parameter extension of (\ref{dem})  is
$$ 
\delta_\mu^{(H,q)}(n):= \sup\set{\delta>0 : \exists\,Q_1,\ldots,Q_n\subset T,\;
J_\mu^{(H,q)}(Q_i)\ge\delta}
$$
where the cubes $Q_i$ are supposed to possess disjoint interiors.
\medskip

Here we shall prove  the following result.
\begin{thm} \label{t:updh}
Let $\mu$ be a measure on a bounded set $T\subset \R^N$ and let $W_H$
be an
$N$--parameter fractional Brownian motion
with Hurst parameter $H$. Let
$q\in[1,\infty)$.
If
$$
 \delta^{(H,q)}_\mu(n)\succeq n^{-1/q-1/a}\,(\log n)^\beta
$$
for certain $a>0$ and $\beta\in\R$, then
$$
 - \log\P\{\norm{W_H}_{L_q(T,\mu)} < \eps\}
 \succeq
 \eps^{-a}\cdot\log(1/\eps)^{a \beta}\;.
$$
\end{thm}
\bigskip
Note that this result does not follow from Theorem
\ref{t:upd} directly, since inequality (\ref{LND})
only holds for $0\le \tau\le \abs{t}$. But we will show that
the proof, based on Theorem \ref{t:upd}, is almost immediate.
\bigskip

We now turn to lower estimates of $\ph_{q,\mu}(\eps)$.
Let $\mu$ be a finite compactly supported Borel measure on $\R^N$.
For a bounded measurable set $A \subset \R^N$
the quantity  $ J_\mu^{(H,q)}(A)$  was introduced in (\ref{Jmu}). Furthermore,
for $H\in (0,1]$ and $q\in [1,\infty)$
 the number $r$ is now (compare with (\ref{r1})) defined by
$$ 
  \frac 1 r := \frac H N +\frac 1 q\; .
$$
Finally, for $n\in  \N$, as in (\ref{sm})  we set
\be
\label{sig}
  \sigma_\mu^{(H,q)}(n)
  :=\inf\set{\left(\sum_{j=1}^n J_\mu^{(H,q)}(A_j)^r\right)^{1/r} :\;
  T\subseteq \bigcup_{j=1}^n A_j}
\ee
where the $A_j$'s are compact subsets of $\R^N$ and $T$ denotes
the support of $\mu$. With this notation we shall prove the following
multi--parameter extension of
(\ref{low}).

\begin{thm}
\label{mup}
Let $X:= (X(t), \, t\in T)$ be a centered Gaussian random field indexed by a
compact set $T\subset\R^N$ and satisfying
$$
\E |X(t)-X(s)|^2\le c\, |t-s|^{2H}\;,\quad
t,s\in T\;,
$$
for some $0<H\le 1$.
If $\mu$ is a finite measure with support in $T$ such that
for certain $q\in[1,\infty)$, $\nu\ge 0$ and $\beta\in \R$
$$
\sigma_\mu^{(H,q)}(n)\preceq n^{-\nu}\, (\log n)^\beta\;,
$$
 then
$$
-\log\P\{ \norm{X}_{L_q(T,\mu)}\le \eps\}
 \preceq  \eps^{-a}\,
\log(1/\eps)^{a\beta}
$$
where $1/a = \nu + H/N$.
\end{thm}
\medskip

\noindent
\textbf{Problem:} Of course, Theorem \ref{mup} applies in particular to $W_H$. Recall that we have
$\e(\abs{W_H(t)-W_H(s)}^2) =\abs{t-s}^{2 H}$, $t,s\in\R^N$. Yet for general measures $\mu$
and $N>1$ we do not know how  the
quantities $\sigma_\mu^{(H,q)}$ and $\delta_\mu^{(H,q)}$ are related
(recall (\ref{q1}) for $N=1$).
Later on we shall prove a relation similar to (\ref{q1}) for a special
class of measures on $\R^N$, the so--called self--similar measures.
But in the general situation the following question remains open. Let $N>1$. Does
as in Theorem \ref{lls} for $N=1$
$$
 \sigma^{(H,q)}_\mu(n)\succeq n^{-\nu}\,(\log n)^\beta
$$
for certain $\nu\ge 0$ and $\beta\in\R$ always imply
$$
 - \log\P\{\norm{W_H}_{L_q(T,\mu)} < \eps\}
 \succeq
 \eps^{-a}\cdot\log(1/\eps)^{a \beta}
$$
with $a$ as in Theorem \ref{mup}?

\bigskip

The rest of the paper is organized as follows.
Section \ref{s:ub} is devoted to the study of upper estimates for small
deviation probabilities, where Theorems \ref{up1}, \ref{t:upd} and \ref{t:updh} are proved.
In Section \ref{s:lb}, we are interested in lower estimates for small deviation probabilities,
and prove Theorem \ref{mup}. Section \ref{s:ss} focuses on the case of self--similar measures.
Finally,  in Section \ref{s:L-infty} we discuss the $L_\infty$--norm.

\section{Upper estimates for small deviation probabilities}
\label{s:ub}
\setcounter{equation}{0}

This section is divided into four distinct parts. The first three parts are devoted to the proof of
Theorems \ref{up1}, \ref{t:upd} and \ref{t:updh}, respectively. The last part contains some
concluding remarks.

\subsection{Proof of Theorem \ref{up1}}
To prove Theorem \ref{up1}, we shall verify the following quite general upper estimate for small
deviation probabilities. As in the formulation of Theorem \ref{up1},
let $X=(X(t), t\in T)$ be a measurable centered Gaussian
process on a metric space $(T,\rho)$, let $\mu$ be a finite Borel measure on $T$
and for disjoint measurable subsets
$A_1,\ldots,A_n$ in $T$ the quantity $V_\mu=V_\mu(A_1,\ldots,A_n)$ is as in (\ref{wmu}).
\begin{prop}
\label{p:ub}
There exist a constant $c_1\in (0,\infty)$ depending only on $q$ and a numerical constant
$c_2\in (0,\infty)$ such that
   $$
   \P\left( \|X\|_{L_q(T,\mu)}^q \le c_1 \, n\, V_\mu^q\right)
   \le \ee^{-c_2\, n}.
   $$
\end{prop}
\pf
For the sake of
clarity, the proof is divided into three distinct steps.

\bigskip

\noindent {\it Step 1. Reduction to independent processes.} We
define the predictions ${\widehat X}_1(t) := 0$, $t\in A_1$, and
${\widehat X}_i(t) := \e\{ X(t) \, | X(s), \, s\in \cup_{k=1}^{i-1}
A_k \}$, $t\in A_i$ (for $2\le i\le n$). The prediction errors are
$$
X_i(t) := X(t) - {\widehat X}_i(t), \qquad t\in A_i, \; 1\le i\le n.
$$

\noindent It is easy to see that $(X_i(t), t\in A_i)_{1\le i\le n}$
are $n$ independent processes: for any $1\le i\le n$, the random
variable $X_i(t)$ is orthogonal to the span of $(X(s), \, s\in
\cup_{k=1}^{i-1} A_k)$ for any $t\in A_i$,
whereas all the random variables $X_j(u)$,
$u\in A_j$, $j<i$, belong to this span.

The main ingredient in Step 1 is the following inequality.

\medskip

\begin{lemma}
 \label{l:anderson}
   For any $\varepsilon>0$,
   $$
   \P\left( \, \sum_{i=1}^n \int_{A_i} |X(t)|^q \d\mu(t) \le
   \eps\right) \le
\P\left( \, \sum_{i=1}^n \int_{A_i}
   |X_i(t)|^q \d\mu(t) \le \eps\right) .
   $$
\end{lemma}

\medskip

\noindent {\it Proof of Lemma \ref{l:anderson}:} There is nothing to
prove if $n=1$. Assume $n> 1$. Let
\begin{eqnarray*}
    {\mathcal F}_{n-1}
 &:=& \sigma\left(  X(s), \, s\in \bigcup_{i=1}^{n-1} A_i \right) ,
    \\
    S_{n-1}
 &:=& \sum_{i=1}^{n-1} \int_{A_i} |X(t)|^q \d\mu(t),
    \\
    U_n
 &:=& \int_{A_n} |X(t)|^q \d\mu(t).
\end{eqnarray*}

\noindent It follows that
\begin{eqnarray}
    \P\left( \, \sum_{i=1}^n \int_{A_i} |X(t)|^q \d\mu(t)
    \le
    \eps \right)
 &=& \P\left( S_{n-1} + U_n \le \eps \right)
    \nonumber
    \\
 &=& \e \left\{ \P\left( S_{n-1} + U_n \le \eps \, \Big|
    \, {\mathcal F}_{n-1} \right) \right\} .
    \label{anderson-split}
\end{eqnarray}

\noindent By definition,
$U_n = \int_{A_n} |X_n(t) + {\widehat X}_n(t) |^q \d\mu(t)
= \| {\widehat X}_n + X_n \|_{L_q(A_n,\mu)}^q$.
Observe that $({\widehat X}_n(t), \, t\in A_n)$ and
$S_{n-1}$ are ${\mathcal F}_{n-1}$-measurable, \ whereas
$(X_n(t), \, t\in A_n)$ is independent of ${\mathcal F}_{n-1}$.
Therefore, by Anderson's inequality (see \cite{And} or \cite{Lif2}),
\begin{eqnarray*}
    \P\left( S_{n-1} + U_n \le \eps \, \Big| \, {\mathcal F}_{n-1}
    \right)
 &\le& \P\left( \| X_n \|_{L_q(A_n, \mu)}^q \le (\eps -
    S_{n-1})_+ \, \Big| \, {\mathcal F}_{n-1} \right)
    \\
 &=& \P\left( \| X_n \|_{L_q(A_n, \mu)}^q + S_{n-1} \le \eps
    \, \Big| \, {\mathcal F}_{n-1} \right) .
\end{eqnarray*}

\noindent Plugging this into (\ref{anderson-split}) yields that
$$
\P\left( \, \sum_{i=1}^n \int_{A_i} |X(t)|^q \d\mu(t) \le
\eps \right) \le \P\left( \| X_n \|_{L_q(A_n,
\mu)}^q + S_{n-1} \le \eps \right) .
$$

\noindent Since the process $(X_n(t), \, t\in A_n)$ and the random
variable $S_{n-1}$ are independent, Lemma \ref{l:anderson} follows by
induction.\hfill$\Box$

\bigskip

\noindent {\it Step 2. Evaluation of independent processes.} In this
step, we even do not use the specific definition of the processes
$X_i(\cdot)$.

\medskip

\begin{lemma}
 \label{l:ub}
   Let $(X_i(t), \, t\in A_i)_{1\le i\le n}$ be independent centered
   Gaussian processes defined on disjoint subsets $(A_i)_{1\le i\le
   n}$ of $T$. Then
   $$
   \P\left( \, \sum_{i=1}^n \int_{A_i} |X_i(t)|^q \d\mu(t)
   \le
   c_1 \, n {\widetilde V_\mu}^q \right) \le \ee^{-c_2\, n},
   $$
   where $c_1$ depends only on $q$, $c_2$ is a numerical constant,
   and
   \be
   \label{tV}
   {\widetilde V_\mu} := \min_{1\le i\le n} \inf_{t\in A_i}
   \left\{ \hbox{\rm Var} (X_i(t)) \right\}^{1/2} \mu(A_i)^{1/q} .
   \ee
\end{lemma}

\medskip

\noindent {\it Proof of Lemma \ref{l:ub}:} Write
$$
Y_i := \int_{A_i} |X_i(t)|^q \d\mu(t), \qquad 1\le i\le n,
$$

\noindent which are independent random variables. We reduce
$\sum_{i=1}^n Y_i$ to a sum of {\it Bernoulli} random variables.
Let
$S_i := Y_i^{1/q}$ and $m_i := \mathrm{median}(S_i)$, $1\le i\le
n$. Consider random variables
$$
B_i := {\bf 1}_{ \{ Y_i \ge m_i^q\} }, \qquad 1\le i\le n .
$$

\noindent Since $m_i^q$ is a median for $Y_i$, we have $\P(B_i=0) =
\P(B_i=1) = 1/2$. In other words, $(B_i, \, 1\le i\le n)$ is a collection
of i.i.d.~Bernoulli random variables.

Since $Y_i \ge m_i^q B_i$, we have, for any $x>0$,
$$
\P\left( \, \sum_{i=1}^n Y_i \le x\right) \le \P\left( \,
\sum_{i=1}^n m_i^q B_i \le x\right)
 \le \P\left( \, \sum_{i=1}^n B_i
\le {x\over \min_{1\le i\le n} m_i^q} \right) .
$$

\noindent In order to evaluate $\min_{1\le i\le n} m_i^q$, we use
the following general result.
\medskip

\begin{fact}
 \label{f:median}
   Let $(X(t), \, t\in T)$ be a Gaussian random process. Assume that\\
   $S:= \sup_{t\in T} |X(t)|<\infty$~a.s. Let $m$ be a median of the
   distribution of $S$. Then
   $$
   m\le \e(S) \le c \, m,
   $$
   where $c:= 1+ \sqrt{2\pi}\,$.
\end{fact}

The first inequality in Fact
\ref{f:median} is in Lifshits~\cite{Lif2}, p.~143, the
second in Ledoux and Talagrand~\cite{LT}, p.~58.
\medskip

Let us complete the proof of Lemma \ref{l:ub}. By Fact \ref{f:median}, we have
$m_i \ge c^{-1} \, \e(S_i) = c^{-1}\, \e(\|X_i\|_{L_q(A_i, \mu)})$. Recall
(Ledoux and Talagrand~\cite{LT}, p.~60) that there
exists a constant $c_q \in (0,\infty)$, depending only on $q$, such
that $\e(\|X_i\|_{L_q(A_i, \mu)}) \ge c_q \{
\e(\|X_i\|_{L_q(A_i, \mu)}^q ) \}^{1/q}$ and that $\e( |X_i(t) |^q)
\ge c_q^q\, \{ \mathrm{Var} (X_i(t))\}^{q/2}$, $t\in A_i$.
Accordingly,
\begin{eqnarray*}
    m_i^q
 &\ge& c^{-q} c_q^q \, \e \left( \|X_i\|_{L_q(A_i, \mu)}^q \right)
    \\
 &=& c^{-q} c_q^q \, \int_{A_i} \e(|X_i(t)|^q) \,
    \d\mu(t)
    \\
 &\ge& c^{-q} c_q^{2q} \, \mu(A_i) \inf_{t\in A_i} \left\{ \mathrm{Var}
    (X_i(t)) \right\}^{q/2} .
\end{eqnarray*}

\noindent Thus $\min_{1\le i\le n} m_i^q \ge c^{-q} c_q^{2q} \,
{\widetilde V_\mu}^q$. It follows that
$$
\P\left( \, \sum_{i=1}^n Y_i \le x\right) \le \P\left( \,
\sum_{i=1}^n B_i \le {x\over c^{-q} c_q^{2q} \,
{\widetilde V_\mu}^q} \right) .
$$

\noindent Taking $x:= {c^{-q} c_q^{2q}\over 3} \,
{\widetilde V_\mu}^q n$, we obtain, by Chernoff's inequality,
$$
\P\left( \, \sum_{i=1}^n Y_i \le {c^{-q} c_q^{2q}\over 3} \,
{\widetilde V_\mu}^q n \right) \le \P\left( \,
\sum_{i=1}^n B_i \le {n\over 3} \right) \le \ee^{-c_2\, n}.
$$

\noindent Lemma \ref{l:ub} is proved, and Step 2
completed.\hfill$\Box$

\bigskip

\noindent {\it Step 3. Final calculations.} We apply the result of
Step 2 to the processes $(X_i(t), \, t\in A_i)$ constructed in Step
1. For any $2\le i\le n$ and any $t\in A_i$, we have
\begin{eqnarray*}
    \mathrm{Var} \left( X_i(t) \right)
 &=& \mathrm{Var}\left( X(t) - {\widehat X}_i(t)  \right)
    \\
 &=& \mathrm{Var} \left[ X(t) - \e \left( X(t) \, \Big| \, X(s) ,
    \, s\in \bigcup_{k=1}^{i-1} A_k \right) \right]
    \\
 &\ge& \mathrm{Var} \left[ X(t) - \e \left( X(t) \, \Big| \, X(s) ,
    \, \rho(s,t) \ge \mathrm{dist} (t, \, \bigcup_{k=1}^{i-1} A_k
    ) \right) \right]
    \\
 &\ge& \mathrm{Var}\left[ X(t) - \e \left( X(t) \, \Big| \, X(s) ,
    \, \rho(s,t) \ge \mathrm{dist} (A_i, \, \bigcup_{k=1}^{i-1} A_k
    ) \right) \right]
    \\
&=& v \left(t,\tau_i\right)^2 ,
\end{eqnarray*}
where $v(t,\cdot)$ is as in (\ref{vt}) and $\tau_i$  as in (\ref{Jbar}).
Therefore, letting $v_i$ be as in (\ref{Jbar}), we get
$$
v_i \le \inf_{t\in A_i} \left\{ \mathrm{Var} (X_i(t))
\right\}^{1/2} \mu(A_i)^{1/q},\qquad
2\le i\le n.
$$
Moreover,  from $X(t)=X_1(t)$, $t\in A_1$, it follows  that
$$
v_1 = \inf_{t\in A_1} \left\{ \mathrm{Var} (X_1(t))
\right\}^{1/2} \mu(A_1)^{1/q}
$$
as well.
Hence, by (\ref{tV}) we get $ V_\mu = V_\mu(A_1,\ldots, A_n) \le {\widetilde V_\mu}$.
It follows from Lemmas \ref{l:anderson} and \ref{l:ub} that
\begin{eqnarray*}
\P\left( \|X\|_{L_q(T,\mu)}^q \le c_1 \, n V_\mu^q \right)
   &\le&
   \P\left( \, \sum_{i=1}^n \int_{A_i} |X(t)|^q\, \d\mu(t)
      \le c_1 \, n V_\mu^q\right)
\cr
    &\le&
    \P\left( \, \sum_{i=1}^n \int_{A_i} |X_i(t)|^q\, \d\mu(t)
      \le c_1 \, n V_\mu^q\right)
\cr
    &\le&
    \P\left( \, \sum_{i=1}^n \int_{A_i} |X_i(t)|^q\, \d\mu(t)
      \le c_1 \, n {\widetilde V_\mu}^q\right)
    \le \ee^{-c_2 \, n}.
\end{eqnarray*}

\noindent This completes Step 3, and thus the proof of Proposition
\ref{p:ub}.\hfill$\Box$
\bigskip

\noindent
\textit{Proof of Theorem \ref{up1}:}
By assumption there is a constant $c>0$ such that
$$
  \delta_n := c\, n^{-1/q-1/a}\,(\log n)^\beta< \delta_\mu(n) \,,\quad n\in\N\;.
$$
Consequently, in view of the definition
of $\delta_\mu(n)$ there exist disjoint measurable subsets $A_1,\ldots, A_n$ in $T$ with
$V_\mu=V_\mu(A_1,\ldots,A_n)\ge \delta_n$. From Proposition \ref{p:ub}, we derive
\be
\label{u1}
\P\{\norm{X}_{L_q(T,\mu)}\le c_1^{1/q}\,n^{1/q}\,\delta_n\}
\le
\P\{\norm{X}_{L_q(T,\mu)}\le c_1^{1/q}\,n^{1/q}\,V_\mu\}\le \mathrm{e}^{-c_2 n}\;.
\ee
Letting $\eps= c_1^{1/q}\,n^{1/q}\,\delta_n = c_1^{1/q}\,c\,n^{-1/a}\,(\log n)^\beta$, it follows that
$c_2\, n\succeq \eps^{-a}\,\log(1/\eps)^{a\beta}$, hence (\ref{u1}) completes the proof of
Theorem \ref{up1}.\hfill$\Box$

\subsection{Proof of Theorem \ref{t:upd}}
\textit{Proof of Theorem \ref{t:upd}:}
This follows from Theorem \ref{up1} and the next
proposition.

\begin{prop}
\label{cubes}
Let $T\subset\R^N$ and let $\mu$ be a finite measure on $T$. Then for $n\in\N$,
$$
\overline \delta_\mu(n)\le 2^{N/q}\,\delta_\mu\big([2^{-N}\,n]\big)
$$
where, as usual, $[x]$ denotes the integer part of a real number $x$.
\end{prop}
\pf
Let $Q_1,\ldots,Q_n$ be arbitrary  cubes in $T$ possessing disjoint interiors.
 Without loss of generality, we may assume that
the diameters of the $Q_i$ are non--increasing.
Set $G := \{ -1, \, 1\}^N$. We cut every cube $Q_i$ into a union
of $2^N$ smaller cubes (by splitting each side into two equal
pieces):
$$
Q_i = \bigcup_{g\in G} Q_i^g.
$$

\noindent For any $i\le n$, let $g(i) \in G$ be such that
$$
\mu(Q_i^{g(i)}) \ge {\mu(Q_i) \over 2^N}
$$
(if the choice is not unique, we choose any one possible value). Let $g\in G$ be such that
\begin{equation}
    \# \left\{ i\in [1,n]\cap \n: \; g(i)= g\right\} \ge {n\over
    2^N}.
    \label{n}
\end{equation}

\noindent We write $I_g := \{ i\in [1,n]\cap \n: \; g(i)= g
\}$, and consider the family of sets $A_i := Q_i^g$, $i\in I_g$. The
following simple geometric lemma provides a lower bound for
$\mathrm{dist} (A_i, \, A_j)$, $i\not=j$.

\medskip

\begin{lemma}
 \label{l:cube}
   Let $Q^\pm := [-1, 1]^N$ and $Q^+ := [0,1]^N$.
   Let $x_1$, $x_2\in \R^N$ and $r_1$, $r_2\in \R_+$ be such that
   the cubes $Q_i := x_i + r_iQ^\pm$, $i=1$ and $2$, are disjoint.
   Then
   $$
   \hbox{\rm dist} (x_1+r_1 Q^+, \, x_2 + r_2 Q^+) \ge \min\{ r_1, \,
   r_2\}.
   $$
\end{lemma}

\medskip

\noindent {\it Proof of Lemma \ref{l:cube}:} For any $x\in \R^N$, we
write $x = (x^{(1)}, \dots, x^{(N)})$. Since the cubes $x_1+r_1
Q^\pm$ and $x_2 + r_2 Q^\pm$ are disjoint, there exists $\ell \in
[1,N] \cap \n$ such that the intervals $[x_1^{(\ell)} - r_1, \,
x_1^{(\ell)} + r_1]$ and $[x_2^{(\ell)} - r_2, \, x_2^{(\ell)} +
r_2]$ are disjoint (otherwise, there would be a point
belonging to both cubes $Q_1$ and $Q_2$). Without loss of generality,
we assume that $x_1^{(\ell)} + r_1 < x_2^{(\ell)} - r_2$. Then, for
any $y_1 \in x_1 + r_1 Q^+$ and $y_2 \in x_2 + r_2 Q^+$, we have
$$
|y_2-y_1| \ge |y_2^{(\ell)} - y_1^{(\ell)}| \ge y_2^{(\ell)} -
y_1^{(\ell)} \ge x_2 ^{(\ell)} - (x_1^{(\ell)} +r_1) \ge r_2 \ge
\min\{ r_1, \, r_2\},
$$

\noindent proving the lemma.\hfill$\Box$

\bigskip

We continue with the proof of Proposition \ref{cubes}. It follows from
Lemma \ref{l:cube} that for any $i>k$ with $i\in I_g$ and $k\in I_g$,
$$
\mathrm{dist} (A_i, \, A_k) \ge N^{-1/2} \min\left\{ \mathrm{diam}
(A_i), \, \mathrm{diam} (A_k) \right\} = N^{-1/2} \mathrm{diam}
(A_i) ,
$$

\noindent (by recalling that the diameters of $Q_i$ are
non--increasing). Let $i_0$ be the minimal element of $I_g$. Then,  for
$i \in I_g,\;i>i_0$,
\begin{eqnarray*}
 v \left(A_i, \mathrm{dist} (A_i, \,
    \bigcup_{k\in I_g,k<i} A_k ) \right)\,\mu(A_i)^{1/q}
 &\ge& v \left( A_i, {\mathrm{diam} (A_i) \over \sqrt{N}} \right)\,\mu(A_i)^{1/q}
    \\
 &\ge& v \left( Q_i, {\mathrm{diam} (Q_i) \over 2\sqrt{N}} \right)\,\mu(Q_i)^{1/q}
     \, 2^{-N/q}
    \\
 &=& 2^{-N/q}\, \overline v(Q_i)\,\mu(Q_i)^{1/q}\\
 &\ge&  2^{-N/q}\,\overline V_\mu(Q_1,\ldots,Q_n)\;.
\end{eqnarray*}
Similarly, using the inequality $\mathrm {Var} X \ge  \mathrm {Var} [ X - \E ( X \, | \, {\cal F} ) ]$
(for any random variable $X$ and any $\sigma$-field ${\cal F}$), we obtain
\begin{eqnarray*}
 \inf_{t\in A_{i_0}} \left(\mathrm {Var} X(t)\right)^{1/2}\cdot \mu(A_{i_0})^{1/q}
&\ge& v \left( A_{i_0}, {\mathrm{diam} (A_{i_0}) \over \sqrt{N}} \right)\,\mu(A_{i_0})^{1/q}
\\
 &\ge& v \left( Q_{i_0}, {\mathrm{diam} (Q_{i_0}) \over 2\sqrt{N}} \right)\,\mu(Q_{i_0})^{1/q} \, 2^{-N/q}
\\
 &=& 2^{-N/q}\, \overline v(Q_{i_0})\,\mu(Q_{i_0})^{1/q}  \ge  2^{-N/q}\,\overline V_\mu(Q_1,\ldots,Q_n)\;.
    \end{eqnarray*}

\noindent
Note that the
cardinality  of $I_g$ which takes the place of the parameter $n$ in $\delta_\mu$ is,
according to (\ref{n}),  not smaller than $2^{-N}n$. Hence, since the cubes $Q_1,\ldots,Q_n$
were chosen arbitrarily in $T$, the proof of Proposition \ref{cubes} follows by the definition of
$\delta_\mu$ and $\overline \delta_\mu$ in (\ref{tilded}) and (\ref{demd}), respectively.
 \hfill$\Box$
\subsection{Proof of Theorem \ref{t:updh}}
\textit{Proof of Theorem \ref{t:updh}: }
Let $T\subset\R^N$ be a bounded set and let $\mu$ be a finite measure on $T$. We
first suppose that $T$ is ``far away from zero", i.e., we assume
\be
\label{distT}
\mathrm{diam}(T)\le \mathrm{dist}(\{0\},T)\;.
\ee
By (\ref{LND}), for any $t\in T$,
$$
v(t,\tau)^2=\mathrm{Var} \left[ W_H(t)- \E\left( W_H(t) \, | \, W_H(s), \,s\in T,\; |s-t|
\ge \tau\right) \right]
\ge c\ \tau^{2H}
$$
for all $\tau\le\mathrm{dist}(\{0\},T)$. Consequently, for any cubes $Q_1,\ldots,Q_n$ in $T$ with
disjoint interiors, we obtain
$$
\overline{V}_\mu (Q_1,\ldots,Q_n)\ge c'\,\min_{1\le i\le n}J_\mu^{(H,q)}(Q_i)\;,
$$
hence $\overline\delta_\mu(n) \ge c'\,\delta_\mu^{(H,q)}(n)$. Theorem \ref{t:updh} follows now from
Theorem \ref{t:upd}.

Next let $T$ be an arbitrary bounded subset of  $\R^N$ and $\mu$ a finite measure on $T$.
We choose an element $t_0\in \R^N$ such that $T_0:= T+t_0$ satisfies (\ref{distT}). By what we have
just proved,
$$
-\log\P\{\norm{W_H}_{L_q(T_0,\mu_0)}<\eps\}\succeq \eps^{-a}\log(1/\eps)^{a\beta}\; ,
$$
where $\mu_0:= \mu* \delta_{\{t_0\}}$ ($\delta_{\{t_0\}}$ being the Dirac measure at $t_0$).
Observe that $\norm{W_H}_{L_q(T_0,\mu_0)} = \{ \int_T\abs{W_H(t+t_0)}^q \d \mu(t) \}^{1/q}$.
Since $\widetilde W_H:=(W_H(t+t_0)-W_H(t_0),\, t\in\R^N)$ is an $N$--parameter fBm as well,
we finally arrive at:
$$
-\log\P\left\{\int_T\abs{\widetilde W_H(t)+W_H(t_0)}^q\d \mu(t)<\eps^q\right\}
\succeq \eps^{-a}\log(1/\eps)^{a\beta}\;.
$$
Theorem \ref{t:updh} follows
from the weak correlation inequality, see \cite{LS}, proof of Theorem 3.7.
\hfill$\Box$

\subsection{Concluding remarks:}
Suppose that $v(t,\tau)\ge c\,\tau^{H}$ and $\mu=\lambda_N$,
the $N$--dimensional Lebesgue measure. Assuming that the interior
of $T$ is non--empty, we easily get
$\delta_\mu(n)\succeq n^{-(1/q+H/N)}$, hence by Theorem \ref{up1},
\[
 - \log\P\{\norm{X}_{L_q(T,\mu)} < \eps\}
 \succeq   \eps^{-N/H}.
\]

Our estimates are suited rather well for stationary fields. For the
non--stationary ones a logarithmic gap may appear. For example, let
$X$ be an $N$--parameter Brownian sheet with covariance
$$\E X(s)X(t) = \prod_{k=1}^N \min\{s_k,t_k\}.
$$
Then we get   $v(t,\tau)\ge c\, \tau^{N/2}$ for all $\tau<\min_{1\le i\le N} t_i$.
By Theorem \ref{t:upd}, for the Lebesgue measure and, say, the $N$--dimensional unit cube $T$,
$$
 - \log\P\{\norm{X}_{L_q(T,\mu)} < \eps\}
 \succeq
 \eps^{-2},
$$
while it is known that in fact
$$
 - \log\P\{\norm{X}_{L_q(T,\mu)} < \eps\}
 \approx
 \eps^{-2}\,\log(1/\eps)^{2N-2} \; .
$$

We also note that cubes in Theorem \ref{t:upd}
can not be replaced by arbitrary closed convex sets. Indeed,
disjoint ``flat" sets are not helpful in this context, as the following example shows.
Define a probability measure $\mu_0$ on $[0,1]$ by
$$ \mu_0= (1-2^{-h})\ \sum_{k=0}^\infty \sum_{i=1}^{2^k}  2^{-k(1+h)}
\delta_{\{i/2^k\}},
$$
where $h>0$ and $\delta_{\{x\}}$ stands for the Dirac mass at point $x$.
Define a measure on the unit square $T$ by $\mu=\mu_0\otimes \lambda_1$.
For a fixed $k$, by taking $A_i=\{i/2^k\}\times [0,1]$, $1\le i\le 2^k$,
we get $n=2^k$ disjoint sets with $\overline V_\mu(A_1,\ldots,A_n)\approx 2^{-k(1+h)/q}$, whatever
the bound for the interpolation error is.
If Theorem \ref{t:upd} were valid in this setting, we would get
$\overline \delta_\mu(n)\succeq n^{-(1+h)/q}$, and
$$
 - \log\P\{\norm{X}_{L_q(T,\mu)} < \eps\}
 \succeq
 \eps^{-q/h},
$$
while it is known, for example, for the 2--parameter Brownian
motion, that in fact
$$
 - \log\P\{\norm{X}_{L_q(T,\mu)} < \eps\}
 \preceq
 \eps^{-4}  \; .
$$
This would lead to a contradiction whenever $q/h>4$.

\section{Lower estimates for small deviation probabilities}
\label{s:lb}

This section is devoted to the study of lower estimates for small deviation probabilities, and is
divided into three distinct parts. In the first part, we present some basic functional analytic tools,
while in the second part, we establish a result for Kolmogorov numbers of operators with values in
$L_q(T,\mu)$. In the third and last part, we prove Theorem \ref{mup}.

\setcounter{equation}{0}

\subsection{Functional analytic tools}

Let $[E,\|\,\cdot\,\|_E]$ and $[F,\|\,\cdot\,\|_F]$ be Banach spaces and let
$u : E\to F$ be a compact operator. There exist several quantities to
measure the degree of compactness of $u$. We shall need two of them, namely,
the sequences $d_n(u)$ and $e_n(u)$ of Kolmogorov and (dyadic) entropy
numbers, respectively. They are defined by
$$
d_n(u):= \inf\{ \|Q_{F_0} u\| : {F_0}\subseteq F,\,{\rm dim} {F_0} < n\}
$$
where for a subspace ${F_0}\subseteq F$ the operator
$Q_{F_0} : F\to F/{F_0}$ denotes the canonical quotient map
from $F$ onto $F/{F_0}$ .
The entropy numbers are given by
$$
e_n(u) := \inf\{\eps>0 : \exists\, y_1,\ldots,y_{2^{n-1}}\in F
\;\mbox{with}\; u(B_E)\subseteq \bigcup_{j=1}^{2^{n-1}}
(y_j+\eps B_F)\}
$$
where $B_E$ and $B_F$ are the closed unit balls of $E$ and
$F$, respectively. We refer to \cite{CS} and \cite{Pis}
for more information about these numbers.

Kolmogorov and entropy numbers are tightly related by the following
result in
\cite{CKP}.
\begin{prop}
\label{CKP}
Let $(b_n)_{n\ge 1}$ be an increasing sequence tending to infinity and
satisfying
$$
\sup_{n\ge 1}\frac{b_{2 n}}{b_n} :=\kappa <\infty\;.
$$
Then there is a constant $c=c(\kappa)>0$ such that for all compact
operators $u$, we have
$$
\sup_{n\ge 1}\, b_n\, e_n(u) \le c\cdot \sup_{n\ge 1}\, b_n\, d_n(u)\;.
$$
\end{prop}
\bigskip

Let $(T,\rho)$ be a compact metric space. Let $C(T)$ denote as usual
the Banach space of continuous functions on $T$ endowed with the norm
$$
\|f\|_\infty:= \sup_{t\in T}|f(t)|\,,\quad f\in C(T)\;.
$$
If $u$ is an operator from
a Banach space $E$ into $C(T)$, it is said
to be $H$--H\"older for some $0<H\le 1$ provided there is a finite constant
$c>0$ such that
\be
\label{Ho}
\abs{(u x)(t_1)-( u x)(t_2)}\le c\cdot \rho(t_1,t_2)^H
\cdot\norm{x}_E
\ee
for all $t_1,t_2\in T$ and $x\in E$. The smallest possible constant $c$
appearing in (\ref{Ho}) is denoted by $\abs{u}_{\rho,H}$
and we write $\abs{u}_H$ whenever the metric $\rho$ is clearly understood.
Basic properties of $H$--H\"older operators may be found in \cite{CS}.

Before stating the basic result about Kolmogorov numbers of H\"older
operators we need some quantity to measure the size of the compact metric space
$(T,\rho)$.
Given $n\in\N$, the  $n$--th entropy number of $T$ (with respect to the metric $\rho$) is defined by
$$
\eps_n(T):=\inf\set{\eps>0 : \exists\; n\;\mbox{$\rho$--balls of radius}
\;\eps\;\mbox{covering}\;T}\;.
$$
\medskip
Now we may formulate Theorem 5.10.1 in \cite{CS} which will be crucial later on.
We state it in the form as we shall use it.

\begin{thm}
\label{cs1}
Let $\H$ be a Hilbert space and let $(T,\rho)$ be a compact metric space such that
\be \label{e1}
   \eps_n(T)\le \kappa\cdot n^{-\nu}\;,\quad n\in\N\,,
\ee
for a certain $\kappa>0$ and $\nu>0$. Then, if $u: \H \to C(T)$ is $H$--H\"older for some $H\in (0,1]$,
then
\be  \label{e2}
   d_n(u)\le c\cdot\max\set{\norm{u},\abs{u}_H}\cdot
   n^{-1/2-H\,\nu}\;,\quad n\in\N\,,
\ee
where $c>0$ depends on $H$, $\nu$ and $\kappa$.
Here, $\norm{u}$ denotes
the usual operator norm of $u$.
\end{thm}

For our purposes it is important to know how the constant $c$ in (\ref{e2})
depends on the number $\kappa$ appearing in (\ref{e1}).
\begin{cor}
\label{cs2}
Under the assumptions of Theorem $\ref{cs1}$, it follows that
$$
   d_n(u)\le c\cdot\max\set{\norm{u},\kappa^H\,\abs{u}_H}\cdot
   n^{-1/2-H\,\nu}\;,\quad n\in\N\,,
$$
with $c>0$ independent of $\kappa$.
\end{cor}
\pf
We set
$$
  \widetilde \rho(t_1,t_2):= \kappa^{-1}\cdot \rho(t_1,t_2)\;,\quad t_1,t_2\in T\,.
$$
If $\widetilde\eps_n(T)$ are the entropy numbers of $T$ with respect to $\widetilde \rho$, then
$\widetilde\eps_n(T)=\kappa^{-1}\cdot \eps_n(T)$, hence, by (\ref{e1}) we have
$\widetilde\eps_n(T)\le n^{-\nu}$, for $n\in\N$. Consequently, an application of Theorem \ref{cs1} yields
\be
\label{e6}
d_n(u)\le c\cdot\max\set{\norm{u},\abs{u}_{\widetilde \rho,H}}
\cdot n^{-1/2-H\,\nu}\;,\quad n\in\N\,,
\ee
where now $c>0$ is independent of $\kappa$. Observe that a
change of the metric does not change
the operator norm of $u$. The proof of the corollary is completed
by (\ref{e6}) and the observation
$\abs{u}_{\widetilde \rho,H} = \kappa^H\cdot \abs{u}_{\rho,H}$.\hfill$\Box$

\subsection{Kolmogorov numbers of operators with values in $L_q(T,\mu)$}

We now state and prove the main result of
this section. Recall that $\sigma_\mu^{(H,q)}(n)$ was defined in (\ref{sig}).
\begin{thm}
\label{maina}
Let $\mu$ be as before a Borel measure on $\R^N$ with compact support $T$
 and let $u$ be an $H$--H\"older operator from a Hilbert space
$\H$ into $C(T)$. Then for all $n,m\in\N$ and $q\in
[1,\infty)$ we have
\be \label{m1}
d_{n+m}\big(u : \H\to
L_q(T,\mu)\big)\le c\cdot\abs{u}_H\cdot
\sigma_\mu^{(H,q)}(m)\cdot n^{-H/N-1/2}\;.
\ee
Here $c>0$ only
depends on $H$, $q$ and $N$. The H\"older norm of $u$ is taken with respect to the
Euclidean distance in $\R^N$.
\end{thm}

\pf Choose arbitrary compact sets $A_1,\ldots,A_m$ covering $T$, the support
of $\mu$. In each $A_j$ we take a fixed element $t_j$, $1\le j\le m$, and
define operators $u_j : \H\to C(A_j)$ via
$$
(u_j h)(t):= (u h)(t) - (u h)(t_j)\,,\quad t\in A_j,\; h\in \H\;.
$$
Thus
\beaa
  \norm{u_j h}_\infty
  &=& \sup_{t\in A_j}\abs{(u_j h)(t)}=
  \sup_{t\in A_j} \abs{(u h)(t)- (u h)(t_j)}\\
  &\le&
  \abs{u}_H\cdot\sup_{t\in A_j} |t-t_j|^H \cdot \norm{h}
  \le
  \diam(A_j)^H \cdot \abs{u}_H\cdot\norm{h}\;,
\eeaa
i.e., the operator norm of $u_j$ can be estimated by
\be \label{m4}
  \norm{u_j}\le \diam(A_j)^H\cdot \abs{u}_H\;.
\ee
Of course,
\be
\label{a1}
|u_j|_H\le |u|_H\;,
\ee
and, moreover, since $A_j\subseteq \R^N$,
\be
\label{a2}
\eps_n(A_j)\le c \cdot {\rm diam} (A_j) \cdot n^{-1/N}
\ee
with some constant $c>0$ depending only on $N$.

Let
\be
\label{a3}
w_j := \mu(A_j)^{1/q}\cdot u_j\;,\quad 1\le j\le m\;.
\ee
An application of Corollary \ref{cs2} with $\nu=1/N$, together with (\ref{m4}),
(\ref{a1}), (\ref{a2}) and (\ref{a3}), yields
\be \label{s1}
   d_n(w_j : \H\to C(A_j))
   \le c\cdot \abs{u}_H\cdot \diam(A_j)^H \cdot\mu(A_j)^{1/q}\cdot  n^{-H/N-1/2}\;,
   \quad n\in \N\;.
\ee
Let $E^q$ be the $\ell_q$--sum of the Banach spaces $C(A_1),\ldots, C(A_m)$, i.e.,
$$
E^q :=\set{(f_j)_{j=1}^m : f_j\in C(A_j)}
$$
and
$$
\|(f_j)_{j=1}^m\|_{E^q}:= \left(\sum_{j=1}^m \|f_j\|_\infty^q\right)^{1/q}\;.
$$
Define $w^q: \H \to E^q$ by
$$
 w^q h := (w_1 h,\ldots,w_m h)\;,
 \quad h\in \H\;.
$$
Proposition 4.2 in \cite{LLS} applies, and (\ref{s1}) leads to
\be \label{uq}
  d_n(w^q)\le c\cdot \abs{u}_H\cdot \left(\sum_{j=1}^m
  \diam(A_j)^{H r}\cdot \mu(A_j)^{r/q}\right)^{1/r}
  \cdot n^{-H/N-1/2}
\ee
where
$$
 1/r= (H/ N + 1/2) - 1/2+ 1/q
         = H/N + 1/q\;.
$$
To complete the proof, set  $B_1=A_1$
and $B_j= A_j\setminus \bigcup_{i=1}^{j-1} A_i$, $2\le j\le m$.
If the operator $\Phi$ from $E^q$ into
$L_q(T,\mu)$ is defined by
$$
  \Phi\big((f_j)_{j=1}^m\big)(t):= \sum_{j=1}^m f_j(t)\cdot
  \frac{\on_{B_j}(t)}{\mu(A_j)^{1/q}}\;,
$$
then $\norm{\Phi}\le 1$ and
\be
\label{uq1}
  \Phi\circ w^q = u - u_0
\ee
where
$$
(u_0 h)(t) = \sum_{j=1}^m (u h)(t_j) \on_{B_j}(t)\,,\quad t\in T\;.
$$
The operator $u_0$
from $\H$ into
$L_q(T,\mu)$ has rank less or equal than $m$. Hence $d_{m+1}(u_0)=0$
and therefore, from algebraic properties of the Kolmogorov numbers and (\ref{uq1})
and (\ref{uq}), it follows that
\beaa
  d_{n+m}(u :\H\to L_q(T,\mu))&\le& d_n(\Phi\circ w^q)+ d_{m+1}(u_0)
\le d_n(w^q)\cdot\norm{\Phi}\\
&\le&
  c\cdot \abs{u}_H\cdot \left(\sum_{j=1}^m \diam(A_j)^{H r}\cdot
  \mu(A_j)^{r/q}\right)^{1/r}
  \cdot n^{-H/N-1/2}\;.
\eeaa
Taking the infimum over all coverings
$A_1,\ldots,A_m$ of $T$ yields (\ref{m1}).\hfill$\Box$

\begin{cor}
\label{entr}
Let $\mu$ be a finite measure on $\R^N$ with support $T$
and suppose that the operator $u$ from the Hilbert space $\H$ into $ C(T)$
is $H$--H\"older for some $H\in (0,1] $.
If
$$
\sigma_\mu^{(H,q)}(n)\le c\cdot n^{-\nu}\cdot (\log n)^\beta
$$
for certain $\nu\ge 0$ and $\beta\in\R$,
then
$$
e_n(u : \H\to L_q (T,\mu)) \le c\cdot \abs{u}_H\cdot n^{-\nu -H/N-1/2}
\cdot (\log n)^\beta\;.
$$
\end{cor}
\pf
Apply Theorem \ref{maina} with $m=n$. The assertion follows from
Proposition \ref{CKP}. \hfill$\Box$

\subsection{Proof of Theorem \ref{mup}}

We start with some quite general remarks about Gaussian processes (cf.~\cite{OW}). Let
$X:= (X(t), \, t\in T)$ be a centered Gaussian process and let us suppose that
$(T,\rho)$ is a compact metric space. Under quite mild conditions, e.g., if
$\rho(t_n, t)\to 0$ in $T$ implies $\e\abs{X(t_n)-X(t)}^2\to 0$,
there are a (separable) Hilbert space $\H$
and an operator $u : \H\to C(T)$ such that
\be
\label{cov}
\e X(t)X(s) = \langle u^* \delta_t, u^* \delta_s\rangle_\H
\ee
where $u^* : C^*(T)\to \H$ denotes the dual operator of $u$ and $\delta_t\in C^*(T)$
is the usual Dirac measure concentrated in $t\in T$. In particular, it follows that
$$
\e\abs{X(t) -X(s)}^2 = \norm{ u^* \delta_t - u^*\delta_s}_\H^2
= \sup_{\norm{h}\le 1}\abs{(u h)(t)-(u h)(s)}^2\;,\quad t,s\in T\;.
$$
Consequently, whenever $u$ and $X$ are related via (\ref{cov}), the operator
$u$ is $H$--H\"older if and only if
\be
\label{Hou}
\left(\e\abs{X(t)-X(s)}^2\right)^{1/2}\le c\cdot \rho(t,s)^H
\ee
for all $t,s\in T$. Moreover, $|u|_H$ coincides with the smallest $c>0$ for
which (\ref{Hou}) holds.
\bigskip

\noindent
{\it Proof of Theorem \ref{mup}:}
We start the proof by recalling a consequence of Theorem 5.1 in \cite{LL}. Suppose
that $u$ and $X$ are related via (\ref{cov}). Then the following are equivalent for any finite Borel
measure $\mu$ on $T$, any $q\in [1,\infty]$,
$a>0$ and $\beta\in\R$:
\begin{enumerate}
\item[(i)]
There is a $c>0$ such that
$$
e_n(u: \H\to L_q(T,\mu))\le c\cdot n^{-1/a -1/2}\,(\log n)^\beta\;.
$$
\item[(ii)]
For some $c>0$ it follows that
$$
-\log\pr{\norm{X}_{L_q(T,\mu)}<\eps} \le c\cdot \eps^a\cdot\log(1/\eps)^{a
\beta}\;.
$$
\end{enumerate}
Taking this into account, Theorem \ref{mup} is a direct consequence of Corollary \ref{entr} and
the above stated equivalence of (\ref{Hou}) with the $H$--H\"older continuity of the corresponding
operator $u$.\hfill$\Box$

\section{Self--similar measures and sets}
\label{s:ss}
\setcounter{equation}{0}

It is a challenging open problem to obtain suitable estimates for $\sigma_\mu^{(H,q)}$ and/or
$\delta_\mu^{(H,q)}$
in the case of arbitrary compactly supported Borel measures $\mu$
on $\R^N$. As already mentioned, we even do not know how  these quantities
are related in the case $N>1$.
Yet if $\mu$ is self--similar, then
suitable estimates for both of these quantities are available.

Let us briefly recall some basic facts about self--similar  measures
which may be found in \cite{Fal} or
\cite{Hut}. An affine mapping $S:\R^N\to \R^N$ is said to
be a contractive similarity provided that
$$
\abs{S(t_1)-S(t_2)} = \lambda \cdot \abs{t_1-t_2}\,,\quad t_1,t_2\in \R^N\,,
$$
with some $\lambda\in (0,1)$. The number $\lambda$ is called the contraction
factor of $S$. Given (contractive) similarities $S_1,\ldots,S_m$ we denote by
$\lambda_1,\ldots,\lambda_m$ their contraction factors.
There exists  a unique compact set $T\subseteq \R^N$
(the self--similar set
generated by the $S_j$'s) such that
$$
  T=\bigcup_{j=1}^m S_j(T)\;.
$$
Let furthermore $\rho_1,\ldots,\rho_m>0$ be weights, i.e., $\sum_{j=1}^m\rho_j=1$.
Then there is a unique Borel probability
measure $\mu$ on $\R^N$  ($\mu$ is called the self--similar measure generated by
the similarities $S_j$ and the weights $\rho_j$) satisfying
$$
 \mu =\sum_{j=1}^m\rho_j\cdot(\mu\circ S_j^{-1})\;.
$$
Note that $T$ and $\mu$ are related via ${\rm supp}(\mu)=T$.

We shall suppose that the similarities satisfy the strong open set condition, i.e.,
we assume that
there exists an open bounded set $\Omega\subseteq \R^N$ with $T\cap\Omega\not=\emptyset$
such that
\be \label{os}
  \bigcup_{j=1}^m S_j(\Omega)\subseteq \Omega\qquad\mbox{and}\qquad S_i(\Omega)\cap
  S_j(\Omega)=\emptyset,\,\;i\not=j\;.
\ee
It is known that then $T\subseteq\overline\Omega$; and since $T\cap\Omega\not=\emptyset$,
we have $\mu(\Omega)>0$,
hence by the results in \cite{LW}, we even have $\mu(\Omega)=1$ and
$\mu(\partial \Omega)=0$. Let us note that
under these assumptions, we have
\be
\label{su1}
\sum_{j=1}^m\lambda_j^N\le 1\;.
\ee

\begin{prop}
\label{cov1}
Let $\mu$ be a self--similar measure generated by similarities $S_j$ with
contraction factors $\lambda_j$ and weights
$\rho_j$, $1\le j\le m$. For $H\in (0,1]$ and $q\in [1,\infty)$,
let $\gamma>0$ be the  unique solution of the equation
\be
\label{ga1}
\sum_{j=1}^m \lambda_j^{H\gamma}\rho_j^{\gamma/q}=1\;.
\ee
Then, under the strong open set condition, we have
$$
\sigma_\mu^{(H,q)}(n) \le c\cdot{\rm diam}(\Omega)^H \cdot n^{-1/\gamma+1/r}
$$
where as before $1/r=H/N+1/q$.
\end{prop}
\pf
By H\"older's inequality and (\ref{su1}), we necessarily have $\gamma\le r$.

We say that $\alpha$ is a word of length $p$ ($p\in\N$) over $\{1,\ldots,m\}$, if
$\alpha=(i_1,\ldots,i_p)$ for certain $1\le i_j\le m$.
For each such word, we define ($\Omega$ being the set appearing in the open set
condition)
\beaa
S_\alpha&:=& S_{i_1}\circ\cdots\circ S_{i_p}\;,\\
\Omega(\alpha)&:=& S_\alpha(\Omega) \; ,\\
\Lambda(\alpha)&:=& (\lambda_{i_1}\cdots\lambda_{i_p})^H
\cdot(\rho_{i_1}\cdots\rho_{i_p})^{1/q}\;.
\eeaa
We need the following estimate.

\begin{lemma}
\label{cover}
For each real number $s>0$, there exist $\ell=\ell(s)$ words
$\alpha_1,\ldots,\alpha_{\ell(s)}$ over  $\{1,\ldots,m\}$ (not necessarily of the same
length) such that the following holds:
\bea
\label{b1}
T
&\subseteq& \bigcup_{i=1}^{\ell(s)}\overline{\Omega(\alpha_i)} \; ,\\
\label{b3}
\max_{1\le i\le \ell(s)} \Lambda(\alpha_i)
&\le& \ee^{-s} \; ,\\
\ell(s)
&\le& c_1\cdot\ee^{\gamma s} \; ,
\label{b2}
\eea
where $\gamma$ was defined by $(\ref{ga1})$.
\end{lemma}

We postpone the proof of Lemma \ref{cover} for a moment, and proceed in the proof
of Proposition \ref{cov1}. Recall that the strong open set condition implies
$\mu(\partial\Omega)=0$; hence for any word $\alpha$, we have
$\mu(\overline{\Omega(\alpha)})=\mu(\Omega(\alpha))$. Accordingly,
\be
\label{b4}
J_\mu^{(H,q)}(\overline{\Omega(\alpha)}) = \Lambda(\alpha)\cdot J_\mu^{(H,q)}(\Omega)
= \Lambda(\alpha)\cdot {\rm diam}(\Omega)^H\;.
\ee
Let $s>0$ be given, and let $\alpha_1,\ldots,\alpha_{\ell(s)}$ be words over  $\{1,\ldots,m\}$
satisfying (\ref{b1}), (\ref{b3}) and (\ref{b2}).
By (\ref{b4}),
$$
\sigma_\mu^{(H,q)}(\ell(s))\le \ell(s)^{1/r} \cdot \ee^{-s} \cdot  {\rm diam}(\Omega)^H\;.
$$
Given $n\in\N$, we define $s>0$ via the equation $c_1\ee^{\gamma s}=n$, where $c_1$ is the constant
in (\ref{b2}). Then $\ell(s)\le n$.
Note that $n\mapsto \sigma_\mu^{(H,q)}(n)$
is non--increasing (since in the definition of $\sigma_\mu^{(H,q)}(n)$, one or several of the
$A_j$ can be empty). Therefore,
$$
\sigma_\mu^{(H,q)}(n)\le \sigma_\mu^{(H,q)}(\ell(s))
\le c\cdot n^{1/r} \cdot n^{-1/\gamma}\cdot {\rm diam}(\Omega)^H
= c\cdot  {\rm diam}(\Omega)^H\cdot n^{-1/\gamma+1/r}
$$
as asserted. \hfill$\Box$
\medskip

\noindent
\textit{Proof of Lemma \ref{cover}:} Let $\Z_+^m$ be the set of vectors $x=(x_1,\ldots,x_m)$
with $ x_j\in \Z$ and $x_j\ge 0$. A path in  $\Z_+^m$ of length $p$ is a sequence $[x^0,\ldots,x^p]$ with
$x^k\in \Z_+^m$. It is said to be admissible provided that $x^0=0$ and for every $k\le p$, there
exists $j_k\le m$ such that $x_{j_k}^k = x_{j_k}^{k-1}+1$ while $x_j^k = x_j^{k-1}$ for all
$j\in \{1,\ldots, m\}\backslash \{j_k\}$. Let
$\mathcal P$ be the set of all admissible paths of any finite length in $\Z_+^m$.

We define a linear function $L: \Z_+^m \to \R$ by
$$
L(x):=  \sum_{j=1}^m d_j\, x_j\;,\quad x=(x_1,\ldots,x_m) \in \Z_+^m\;,
$$
with $d_j:= -\log(\lambda_j^H\cdot \rho_j^{1/q})$, $1\le j\le m$.
Note that $L$ takes values in $[0,\infty)$, and increases along an
admissible path.
For $s>0$,
let $\mathcal P_s$ be the set of those paths  $[x^0,\ldots,x^p]$ in $\mathcal P$
for which $L(x^{p-1})<s\le L(x^p)$. It was shown in \cite{NS} that
$$
\#(\mathcal P_s) \le c_1\cdot \ee^{\gamma s}
$$
for a certain constant $c_1>0$ and with $\gamma$ satisfying
\begin{equation}
\label{ga2}
1= \ee^{-d_1 \gamma}+\cdots+ \ee^{-d_m \gamma}\;.
\end{equation}
In view of the definition of $d_j$, the number $\gamma$ in (\ref{ga2}) coincides with the one defined in
(\ref{ga1}).
To each path $[x^0,\ldots,x^p]$ in $\Z_+^m$, we  assign a word $\alpha=(j_1,\ldots,j_p)$
as follows: $j_k$ is such that $x_{j_k}^k = x_{j_k}^{k-1}+1$. In this way,
we obtain a one--to--one correspondence between $\mathcal P$ and the set of finite words
over $\{1,\ldots,m\}$. Moreover, a path belongs to $\mathcal P_s$ if and only if for the
corresponding word $\alpha$ we have $\Lambda(\alpha)\le \ee^{-s}<\Lambda(\bar \alpha)$
with $\bar\alpha= (j_1,\ldots,j_{p-1})$. We
hereby set $\Lambda(\bar\alpha)=1$ provided $\bar\alpha$ is the empty word.
If we enumerate all words $\alpha$ corresponding to paths in $\mathcal P_s$,
we get $\alpha_1,\ldots,\alpha_{\ell(s)}$ with $\ell(s)\le c_1\cdot \ee^{\gamma s}$
and $\Lambda(\alpha_i)\le \ee^{-s}$, $1\le i\le \ell(s)$. Moreover, since
$T\subseteq \overline\Omega$, it follows that
$$
T= \bigcup_{i=1}^{\ell(s)}S_{\alpha_i}(T)
\subseteq \bigcup_{i=1}^{\ell(s)} S_{\alpha_i}(\overline \Omega)
= \bigcup_{i=1}^{\ell(s)}\, \overline{\Omega(\alpha_i)}
$$
which completes the proof. \hfill$\Box$
\bigskip

As a consequence of Proposition \ref{cov1} and Theorem \ref{maina}, we get the following.

\begin{thm}
Let $\mu$ be as in Proposition $\ref{cov1}$ with support $T\subset \R^N$.
If $u$ is an
$H$--H\"older operator from a Hilbert space $\H$ into $C(T)$, then
$$
  d_n(u : \H\to L_q(T,\mu))\le c\cdot \diam(\Omega)^H\cdot n^{-1/\gamma+1/q-1/2}\; ,
$$
where $\gamma$ is defined by equation $(\ref{ga1})$. The entropy numbers of $u$ can be estimated in
the same way.
\end{thm}
\bigskip

\begin{cor}
\label{corsu}
Let $\mu$ be a self--similar measure with support $T\subseteq\R^N$ and let
$X:= (X(t), \, t\in T)$ be a centered Gaussian process satisfying
$$
\E|X(t)-X(s)|^2  \le c\cdot |t-s|^{2H} \; , \quad t,s\in T\,,
$$
for some $H\in(0,1]$. Then
$$
-\log \pr { \norm{X}_{L_q(T,\mu)} <\eps }
 \preceq  \eps^{-\gamma q/(q-\gamma)}
$$
with $\gamma$ defined by $(\ref{ga1})$.
\end{cor}
\bigskip

Our next aim is to find suitable lower estimates for $\delta_\mu^{(H,q)}(n)$ in
the case of self--similar measures $\mu$ on $\R^N$. To this end, let $\Omega$ be
the open bounded set satisfying (\ref{os}), and consider
$$
\delta_\mu^{(H,q)}(n):= \sup\set{\delta>0 : \exists\,Q_1,\ldots,Q_n\subset \Omega,\;
J_\mu^{(H,q)}(Q_i)\ge\delta}
$$
where the cubes $Q_i$ are supposed to possess disjoint interiors.
\begin{prop}
\label{lowsf}
Let $\gamma>0$ be as in $(\ref{ga1})$. Then
$$
\delta_\mu^{(H,q)}(n) \succeq n^{-1/\gamma} \; .
$$
\end{prop}
\pf
For an open subset $G\subseteq \R^N$ and $\delta>0$, let $M_\mu=M_\mu(\delta,G)$ be
the maximal number of cubes $Q_1,\ldots,Q_{M_\mu}$ in $G$ with disjoint interiors and with
$J_\mu^{(H,q)}(Q_j)\ge\delta$, $1\le j\le M_\mu$.
For $\Omega$ as above, define open sets $\Omega_i$ and measures $\mu_i$
on $\Omega_i$ by
$$
\Omega_i:=S_i(\Omega)\quad\mbox{ and}\quad  \mu_i := \rho_i\cdot(\mu\circ S_i^{-1})\;,\quad
1\le i\le m\;.
$$
{F}rom (\ref{os}), we derive
\be
\label{mm1}
 M_\mu(\delta,\Omega)\ge \sum_{i=1}^m M_{\mu_i}(\delta,\Omega_i)\;.
\ee
If $Q\subseteq \Omega_i$ is a cube, then by self--similarity, we have
$J_{\mu_i}^{(H,q)}(Q)= \lambda_i^H \, \rho_i^{1/q}\,J_\mu^{(H,q)}(Q)$,
hence $M_{\mu_i}(\delta,\Omega_i) = M_\mu(\beta_i\,\delta,\Omega)
$ with $\beta_i:= \lambda_i^{-H}\rho_i^{-1/q}$. By (\ref{mm1}),
$M_\mu(\delta,\Omega) \ge \sum_{i=1}^m  M_{\mu}(\beta_i\,\delta,\Omega)$ for all $\delta>0$.
Applying Lemma 5.1 (first part) in \cite{LLS} yields, for any $\delta_0>0$,
\be
\label{mm2a}
\inf_{\delta\le\delta_0}\delta^\gamma\,M_\mu(\delta,\Omega)
\ge c\, \delta_0^\gamma\,M_\mu(\delta_0,\Omega)
\ee
where $\gamma$ is defined as the unique solution of $\sum_{i=1}^m \beta_i^{-\gamma}=1$
(thus as in (\ref{ga1})) and the constant
$c>0$ only depends
on $\beta_1,\ldots,\beta_m$.

Since $\Omega$ is open and $\mu(\Omega)=1$, there exists a non--empty cube $Q_0\subseteq\Omega$
such that $\delta_0:=J_\mu^{(H,q)}(Q_0)>0$. Thus $M_\mu(\delta,\Omega)\ge 1$ for $\delta\le \delta_0$.
In view of (\ref{mm2a}), we have
$$
M_\mu(\delta,\Omega)\ge c_0\, \delta^{-\gamma}
$$
whenever $0<\delta\le\delta_0$. Take $\delta:=c_0^{-1/\gamma} n^{-1/\gamma}$. We have proved
that when $n$ is sufficiently large, there exist $n$ cubes
$Q_1,\ldots,Q_n$ in $\Omega$ possessing disjoint interiors such that $J_\mu^{(H,q)}(Q_i)\ge
\delta$. This completes the proof.
\hfill$\Box$
\bigskip

\remark
Propositions \ref{cov1} and \ref{lowsf} imply that
$$
n^{1/r}\delta_\mu^{(H,q)}(n)\succeq \sigma_\mu^{(H,q)}(n)
$$
for self--similar measures $\mu$.
\medskip

As a consequence of Theorem \ref{t:updh} and Proposition \ref{lowsf} we get the following.
\begin{cor}
\label{corsf}
 If $\mu$ is self--similar as before (in particular, the strong
open set condition is assumed), then
$$
-\log \pr { \norm{W_H}_{L_q(\R^N,\mu)} <\eps }
 \succeq  \eps^{-\gamma q/(q-\gamma)}
$$
with $\gamma$ defined by $(\ref{ga1})$.
\end{cor}

Combining Corollaries \ref{corsu} and \ref{corsf} finally gives us the following.
\begin{thm}
\label{thmsf}
 If
$\mu$ is a self--similar measure on $\R^N$ as before, then
$$
-\log \pr { \norm{W_H}_{L_q(\R^N,\mu)} <\eps }
 \approx  \eps^{-\gamma q/(q-\gamma)}\;.
$$
\end{thm}

\noindent
\textbf{Example:} Suppose that the weights and the contraction factors in the construction
of $\mu$ are related via
$$
\lambda_i = \rho_i^s\;,\quad 1\le i\le m\,,
$$
for some $s>0$.   Then it follows that $\gamma=(s H +1/q)^{-1}$, hence
\be
\label{mm3}
-\log \pr { \norm{W_H}_{L_q(\R^N,\mu)} <\eps }
 \approx  \eps^{-1/(s H)}
\ee
for these special weights.

Of special interest is the case $s=1/D$ where $D$ denotes the similarity dimension of the
self--similar set $T$, i.e.,
\be
\label{mm3a}
\sum_{i=1}^m \lambda_i^D = 1\;.
\ee
Then (cf.~\cite{Hut}) $\mu$ is the (normalized) Hausdorff measure on $T$ and $D$ its Hausdorff
dimension. Thus (\ref{mm3}) becomes
\be
\label{mm4}
-\log \pr { \norm{W_H}_{L_q(\R^N,\mu)} <\eps }
 \approx  \eps^{-D/H}
\ee
for Hausdorff measures $\mu$, as claimed in Theorem \ref{th_exammple}.
\bigskip

\remark
For the Lebesgue measure on $[0,1]^N$ and $q=2$,   a more precise asymptotic for (\ref{mm4}) was recently
evaluated in \cite{NN} by means of Hilbert space methods.

\section{$L_\infty$-norm}
\label{s:L-infty}
\setcounter{equation}{0}
In the case $q=\infty$, the natural setting of our problem is as follows: given
a metric space $(T,\rho)$ and a centered Gaussian process $X:= (X(t), \, t\in T)$, evaluate
\begin{equation} \label{task}
\P \left(
\, \sup_{t\in T} |X(t)| \le \eps \right), \qquad \eps\to 0.
\end{equation}
There is no reasonable place for a measure $\mu$ in this problem.

The main tool for working with (\ref{task}) is provided by packing
and entropy cardinalities defined as follows:
\[
 M(\eps,T) := \max\set{n :
    \exists\   t_1,\ldots,t_n\; \mbox{in} \; T\ \mbox{such that} \
    \rho(t_i,t_j)>\eps, \ i\not=j } ,
\]
\[
 N(\eps,T) := \min\set{n :
    \exists\   t_1,\ldots,t_n\; \mbox{in} \; T\ \mbox{such that} \
   T\subset \cup_{j=1}^n B(t_j,\eps) },
\]
where $B(t,\eps)$ denotes the ball with center $t$ and radius $\eps$.
Recall that the asymptotic behavior of $M(\cdot,T)$ and $N(\cdot,T)$
at zero is essentially the same, since
\[
N(\eps,T) \le M(\eps,T) \le N(\frac \eps 2,T).
\]

In order to establish an upper bound for the small deviation probability in (\ref{task}), let us use
$v(t,\tau)$ and $v(A,\tau)$ as defined in (\ref{vt}) and (\ref{vt2}), respectively.
Let $\psi_1(\tau)= v(T,\tau)$. Denote $\psi_1^{-1}(\cdot)$ the inverse function of $\psi_1$.
Then a trivial argument based on the non-determinism
property (see Proposition 3 in \cite{Lif2}, pp.~20--21) yields
\[
\P \left(
\,\sup_{t\in T} |X(t)| \le \eps \right) \le
\exp\left( - c_1 M\left(\psi_1^{-1}(\eps),T\right) \right),
\]
with some numerical constant $c_1>0$.

In order to establish a lower bound for the small deviation probability, we assume that a
H\"older--type condition holds:
\[
\E \left( |X(t)-X(s)|^2 \right) \le \psi_2(\rho(s,t))^{2}\;,\quad t,s\in T.
\]
Then, under minimal regularity assumptions on $\psi_2(\cdot)$
and $N(\cdot,T)$, such as
\be
\label{inf0}
c_2 \, N\left(\psi_2^{-1}(\eps),T\right)
\le N\left(\psi_2^{-1}(\eps/2),T\right)
\le c_3 \ N\left(\psi_2^{-1}(\eps),T\right),
\ee
with some $c_2$ and $c_3>1$,
Talagrand's lower bound (see the original result in \cite{Tal1} and a
better exposition in \cite{Led}) applies and we have
\be
\label{inf1a}
\P \left(
\,\sup_{s,t\in T} |X(t)-X(s)| \le \eps \right) \ge
\exp\left( - c_4 N\left(\psi_2^{-1}(\eps),T\right) \right),
\ee
with a numerical constant $c_4>0$. Notice that, if necessary, the function
$N(\psi_2(\cdot),T)$ can be replaced by any majorant in the
regularity condition.

If we assume that there exists a non--decreasing function $\psi$ regularly varying
at zero such that
\be \label{reg}
\psi_1\approx\psi_2\approx \psi,
\ee
then, by combining the upper and lower estimates, we get
\be
\label{inf1}
-\log \P \left(
\,\sup_{t\in T} |X(t)| \le \eps \right) \approx
N\left(\psi^{-1}(\eps),T\right).
\ee
For the one--parameter version of this result related  to Riemann--Liouville
processes and fractional Brownian motions, we refer to \cite{LJAT}.
\bigskip

For self--similar sets $T\subseteq\R^N$, (\ref{inf1}) bears a particularly simple form, stated as follows.

\begin{prop}
\label{infsf}
Let $T\subseteq\R^N$ be a compact self--similar set such that
the open set condition holds. If a Gaussian process
$X$ on $\R^N$ satisfies $(\ref{reg})$, then
\be
\label{inf2}
-\log \P \left(
\,\sup_{t\in T} |X(t)| \le \eps \right) \approx
(\psi^{-1}(\eps))^{-D}
\ee
where the constant $D>0$ is defined by equation
$(\ref{mm3a})$
and coincides with the Hausdorff dimension of $T$.
\end{prop}
\pf
By Theorem 1 in \cite{La}, we have
$$
N(\eps,T)\approx \eps^{-D}\;,
$$
and since $\psi$ is regularly varying at zero, so
is $\psi^{-1}$, and (\ref{inf0}) is satisfied. Thus (\ref{inf2})
follows from (\ref{inf1}).\hfill$\Box$
\medskip

Arguing as in the proof of Theorem \ref{mup}, by
(\ref{LND}) the preceding proposition applies to $W_H$ and
$\psi(\tau)=\tau^H$. Consequently, Proposition \ref{infsf} leads to
the following.
\begin{cor}
\label{inft}
Let $T\subset\R^N$ be self--similar such that the open set condition holds.
Then it follows that
$$
-\log \P \left(
\,\sup_{t\in T} |W_H(t)| \le \eps \right) \approx
\eps^{-D/H}
$$
where as before $D>0$ denotes the Hausdorff dimension of $T$.
\end{cor}
\bigskip

\remark
In the case $T=[0,1]^N$ Corollary \ref{inft}
recovers Theorem 1.2 of Shao and Wang~\cite{SW} asserting
$$
\exp\left( - {c_1 \over \eps^{N/H}} \right) \le \P \left(
\,\sup_{t\in [0,1]^N} |W_H(t)| \le \eps \right)  \le
\exp\left( - {c_2 \over \eps^{N/H}} \right).
$$

We also mention interesting small deviation bounds
in the $\sup$--norm for {\it stationary} random fields in \cite{LifTsy} and \cite{Tsy}.
\bigskip

One question still has to be answered when comparing Theorem \ref{mup} with (\ref{inf1a}),
namely, whether or not (\ref{inf1a}) may be viewed (if  $\psi_2(\lambda)=c\,\lambda^H$)
as an extension of
Theorem \ref{mup} to the limit case $q=\infty$. The answer is affirmative.
If $T\subseteq\R^N$, $H>0$, then a
natural generalization of (\ref{sig}) to $q=\infty$ is as follows:
$$
\sigma(n)= \sigma^{(H,\infty)}(n,T):= \inf\set{\left(\sum_{j=1}^n {\rm diam}(A_j)^N
\right)^{H/N} : T\subseteq \bigcup_{j=1}^n A_j} \; .
$$
At a first glance, it is not clear how this quantity is related to $N(\eps,T)$ or
$M(\eps,T)$. We will find a connection in terms of the inverse of
$M(\,\cdot\,,T)$, i.e., in terms of the inner entropy numbers $\delta_n$ of $T$ defined by
\beaa
\delta_n &=& \delta_n(T) := \sup\set{\delta> 0 : \exists\; t_1,\ldots,t_n\in T\;, \abs{t_i-t_j}>\delta\;,
1\le i<j\le n}\\
&=& \inf\set{\delta>0 : M(\delta,T)\le n} \, .
\eeaa
The following proposition relates the sequences $(\sigma(n))$ and $(\delta_{n})$.
\begin{prop}
\label{entrop}
There are positive constants $c_1,c_2$ and an integer $\kappa$ depending only on $H$ and $N$ such
that for $T\subseteq \R^N$,
\be
\label{delkap}
c_1\cdot n^{H/N}\cdot\delta_{\kappa n}^H \le \sigma(n)\le c_2\cdot n^{H/N}\cdot \delta_n^H\;.
\ee
\end{prop}
\pf
Let $\kappa$ be any fixed integer with $\kappa>2^N$ and choose $t_1,\ldots,t_{\kappa n}\in T$
such that $|t_i-t_j| \ge \delta_{\kappa  n}$ for $i\not= j$. Then the open balls
 $B(t_i,\delta_{\kappa n}/2)$
are disjoint. If $A_1,\ldots,A_n$ is any covering of $T$ by compact sets, then for each $i\le \kappa n$
there is a $j\le n$ such that
$$
B(t_i,\delta_{\kappa n}/2) \subseteq A_j + B(0,\delta_{\kappa n}/2)\;.
$$
Let $V_N$ be the volume of the $N$--dimensional Euclidean unit ball. Then
\beaa
\kappa n\cdot V_N\cdot(\delta_{\kappa n}/2)^N &=& \sum_{i=1}^{\kappa n}
{\rm vol}_N(B(t_i,\delta_{\kappa n}/2) )
\\
&=&  {\rm vol}_N\left(\bigcup_{i=1}^{\kappa n} B(t_i,\delta_{\kappa n}/2)\right)
\\
&\le& {\rm vol}_N\left(\bigcup_{j=1}^n  \left(A_j + B(0,\delta_{\kappa n}/2)\right)\right)
\\
&\le& V_N\cdot \sum_{j=1}^n \left({\rm diam}(A_j)+ \delta_{\kappa n}/2\right)^N .
\eeaa

\noindent By means of the elementary inequality $(a+b)^N \le 2^N (a^N + b^N)$ (for $a\ge 0$
and $b\ge 0$), this yields
\begin{eqnarray*}
    \kappa n\cdot V_N\cdot(\delta_{\kappa n}/2)^N
 &\le&2^N \cdot V_N \cdot  \sum_{j=1}^n \left[
    {\rm diam}(A_j)^N +(\delta_{\kappa n}/2)^N
    \right]
    \\
 &=& 2^N \cdot V_N\cdot\sum_{j=1}^n
    {\rm diam}(A_j)^N + V_N\cdot n\cdot
    \delta_{\kappa n}^N\;.
\end{eqnarray*}

\noindent Consequently,
$$
n\cdot\left(\frac{\kappa}{2^N}-1\right)\cdot 2^{-N}\cdot\delta_{\kappa n}^N
\le  \sum_{j=1}^n {\rm diam}(A_j)^N\;.
$$
This being true for all compact coverings $A_1,\ldots,A_n$ of $T$, the first inequality in
(\ref{delkap}) follows.

To prove the second inequality in (\ref{delkap}), we take any $\delta>\delta_n$ and a
maximal number $m$ such that there exist
$t_1,\ldots,t_m\in T$ with $\abs{t_i-t_j}>\delta$, $1\le i\not= j\le m$. From
$\delta>\delta_n$ necessarily follows $m\le n$.
Moreover, by the maximality of $m$ we have $T\subseteq
\bigcup_{j=1}^m \overline{B(\delta,t_j)}$. Note that $\sigma$ is decreasing, thus this
implies
$$
\sigma(n)\le \sigma(m)\le   m^{N/H}\cdot (2\,\delta)^H\le  n^{N/H}\cdot (2\,\delta)^H\;.
$$
Since $\delta>\delta_n$ was chosen arbitrarily, this
completes the proof of (\ref{delkap}).\hfill$\Box$
\bigskip

\bigskip

\noindent {\Large\bf Acknowledgements}
\bigskip

\noindent
The work of the first named author was partially supported
by grants RFBR 05-01-00911, RFBR/DFG 04-01-04000
and INTAS 03-51-5018. The second and third named authors were partially supported by
the program PROCOPE 2005/D/04/27467.
\bigskip

\

\bibliographystyle{amsplain}

\begin{thebibliography}{10}



{\baselineskip=12pt



\bibitem{And}
 T. W. Anderson, \textit{The integral of symmetric
 unimodular functions over a symmetric convex set and some
 probability inequalities.} Proc. Amer. Math. Soc.
 \textbf{6} (1955), 170--176.

\bibitem{CKP}
 B. Carl, I. Kyrezi and A. Pajor,
 \textit{Metric entropy of convex hulls in Banach spaces.}
 J. London Math. Soc. \textbf{60} (1999), 871--896.

\bibitem{CS}

 B. Carl and I. Stephani,
 \textit{Entropy, Compactness and Approximation of Operators.}
 Cambridge Univ. Press, Cambridge, 1990.

\bibitem{Fal}
K. J. Falconer,
\textit{The Geometry of Fractal Sets.}
 Cambridge Univ. Press, Cambridge, 1985.

\bibitem{Hut}
  J. E. Hutchinson,
  \textit{Fractals and self-similarity.}
  Indiana Univ. Math. J. \textbf{30} (1981), 713--747.

\bibitem{La}
S. P. Lalley,
\textit{The packing and covering functions of some self--similar
fractals.}
Indiana Univ. Math. J. \textbf{37} (1988), 699--709.

\bibitem{LW} K. S. Lau and J. Wang,
\textit{ Mean quadratic variations and Fourier asymptotics of
self--similar measures.}
Monat. Math. \textbf{115} (1993), 99--132.

\bibitem{Led}
M. Ledoux, \textit{Isoperimetry and Gaussian analysis.} Lectures on Probability
  Theory and Statistics, Lecture Notes in Math., vol. \textbf{1648}, Springer, 1996,
  pp.~165--294.

\bibitem{LT}
 M. Ledoux, and M. Talagrand,
 \textit{Probability in Banach Spaces.} Springer, Berlin, 1991.

\bibitem{Li1}
 W. V.  Li,
 \textit{Small ball estimates for Gaussian
 Markov processes under $L_p$--norm.}
 Stoch. Proc. Appl. \textbf{92} (2001), 87--102.

\bibitem{LL}
 W. V. Li and W. Linde,
 \textit{Approximation, metric entropy and small ball
 estimates for Gaussian measures.} Ann. Probab. \textbf{27} (1999),
 1556--1578.

\bibitem{LS}
 W. V. Li and Q.-M. Shao,
 \textit{Gaussian processes: inequalities, small ball probabilities and
 applications.} In: Shanbhag, D. N. et al. (eds.),
 Stochastic processes: Theory and methods.
 Handb. Statist. \textbf{19} (2001), 533--597.
 Elsevier, Amsterdam.

\bibitem{Lif2}
 M. A. Lifshits,
 \textit{Gaussian Random Functions.} Kluwer, Dordrecht, 1995.

\bibitem{Lif}
 M. A. Lifshits,
 \textit{Asymptotic behavior of small ball probabilities.}
 In: Probab. Theory and Math. Statist. Proc. VII International
 Vilnius Conference, pp.~453--468.
 VSP/TEV, Vilnius, 1999.

\bibitem{LLMem}
 M. A. Lifshits  and W. Linde,
 \textit{Approximation and entropy numbers of Volterra operators
with application to Brownian motion.}
 Memoirs Amer. Math. Soc. \textbf{745} (2002),  1--87.

\bibitem{LLTAMS}
 M. A. Lifshits and W. Linde,
 \textit{Small deviations of weighted fractional processes
 and average non--linear approximation.}
 Trans. Amer. Math. Soc.  \textbf{357} (2005),  2059--2079.

\bibitem{LLS}
 M. A. Lifshits, W. Linde and Z. Shi,
 \textit{Small deviations for  Riemann-Liouville
 processes in $L_q$-spaces with respect to fractal measures.}
Proc. London Math. Soc. \textbf{92} (2006), 224--250.

\bibitem{LSi}
 M. A. Lifshits and T. Simon,
 \textit{Small ball probabilities for stable Riemann-Liouville
 processes.}
 Ann. Inst. H. Poincar\'e \textbf{41} (2005), 725--752.

\bibitem{LifTsy}
M. A. Lifshits and B. S. Tsyrelson, \emph{Small ball deviations of {G}aussian
  fields.} Theor. Probab. Appl. \textbf{31} (1986), 557--558.

\bibitem{LJAT}
 W. Linde,
 \textit{Kolmogorov numbers of Riemann--Liouville operators over
 small sets and applications to Gaussian processes.}
 J. Appr. Theory \textbf{128} (2004), 207--233.

 \bibitem{OW}
 W. Linde,
 \textit{Small ball problems and compactness of operators.}
In: Mathematisches Forschungsinstitut Oberwolfach
Report No. 44/2003. Mini-Workshop: Small Deviation Problems for
Stochastic Processes and Related Topics.
\ {\tt http://www.mfo.de/programme/schedule/2003/42c/Report44\_2003.pdf}


\bibitem{Ma}
P. Mattila,
\textit{Geometry of Sets and Measures in Euclidean Spaces.}
Cambridge Univ. Press, Cambridge, 1995.

\bibitem{NS}
K. Naimark and M. Solomyak,
\textit{The eigenvalue behaviour for the boundary value problems
related to self-similar measures on $R^d$.}
Math. Research Letters \textbf{2} (1995), 279--298.

\bibitem{Na}
A. I. Nazarov,
\textit{Logarithmic asymptotics of small deviations for
 some Gaussian processes in the $L_2$--norm
with respect to a self-similar measure.} Zap. Nauchn.
Sem. S.-Petersburg. Otdel. Mat. Inst. Steklov. (POMI)
\textbf{311}  (2004), 190--213.

\bibitem{NN}
A. I. Nazarov and Y. Y. Nikitin,
\textit{Logarithmic {$L_2$}-small ball asymptotics for some fractional
Gaussian processes.} Theor. Probab. Appl.
\textbf{49} (2004), 695-711.

\bibitem{Ol}
L. Olsen,
\textit{A multifractal formalism.}
Adv. Math. \textbf{116} (1995), 82--196.

\bibitem{Pe}
Ya. B. Pesin,  \textit {Dimension theory in dynamical systems. Contemporary views
and applications.}  Chicago Lectures in Mathematics. University of Chicago
Press, Chicago, 1997.

\bibitem{Pis}
G. Pisier,
\textit{The Volume of Convex Bodies and Banach Space Geometry.}
Cambridge Univ. Press, Cambridge, 1989.

\bibitem{pitt}
   L. D. Pitt,
    \textit{Local times for Gaussian vector fields.}
    Indiana Univ. Math. J. \textbf{27} (1978), 309--330.
\bibitem{SW}
    Q.-M. Shao  and D. Wang,
    \textit{Small ball probabilities of Gaussian fields.}
    Probab. Theory Related Fields \textbf{102} (1995), 511--517.

\bibitem{Tal1}
M. Talagrand, \textit{New Gaussian estimates for enlarged balls.} Geom.
  Funct. Anal. \textbf{3} (1993), 502--526.

\bibitem{Tsy}
B. S. Tsyrelson, \textit{On stationary Gaussian processes with a finite
  correlation function} J. Math. Sci. \textbf{68} (1994), 597--604.

\bibitem{xiao}
   Y.  Xiao,
    \textit{H\"older conditions for the local times and the Hausdorff
    measure of the level sets of Gaussian random fields.}
    Probab. Theory Related Fields \textbf{109} (1997), 129--157.


}
\end{thebibliography}

\bigskip
\bigskip

{\footnotesize
\baselineskip=12pt
\begin{center}
\noindent
\parbox[t]{5cm}
{Mikhail Lifshits\\
St.Petersburg State University\\
198504 Stary Peterhof\\
Dept of Mathematics \\
and Mechanics\\
Bibliotechnaya pl., 2\\
Russia\\
{\tt lifts@mail.rcom.ru}}
\parbox[t]{5cm}
{Werner Linde\\
FSU Jena \\
Institut f\"ur Stochastik\\
Ernst--Abbe--Platz 2\\
07743 Jena\\
Germany\\
{\tt lindew@minet.uni-jena.de}}
\parbox[t]{5cm}
{Zhan Shi\\Laboratoire de Probabilit\'es\\ et Mod\`eles
Al\'eatoires\\
Universit\'e Paris VI\\
4 place Jussieu\\
F-75252 Paris Cedex 05\\
France\\
{\tt zhan@proba.jussieu.fr}}
\end{center}
}
\end{document}